\documentclass[12pt]{article}%
\usepackage{amssymb}
\usepackage{amsfonts}
\usepackage{amsmath}
\usepackage{graphicx}%
\providecommand{\U}[1]{\protect\rule{.1in}{.1in}}
\newtheorem{theorem}{Theorem}[section]

\newtheorem{conjecture}[theorem]{Conjecture}

\newtheorem{Fact}[theorem]{Fact}

\newtheorem{problem}[theorem]{Problem}
\newtheorem{proposition}[theorem]{Proposition}

\providecommand{\boksie}{\ensuremath{\mathbin{\raisebox{0.3mm}{$\scriptstyle\square$}}}}
\begin{document}

\title{\textbf{Protecting a Graph with Mobile Guards}}
\author{William F. Klostermeyer\\School of Computing\\University of North Florida\\Jacksonville, FL 32224-2669\\{\small wkloster@unf.edu}
\and C. M. Mynhardt\thanks{Supported by the Natural Sciences and Engineering
Research Council of Canada.}\\Department of Mathematics and Statistics\\University of Victoria, P.O. Box 3060 STN CSC\\Victoria, BC, \textsc{Canada} V8W 3R4\\{\small mynhardt@math.uvic.ca}}
\maketitle

\begin{abstract}
Mobile guards on the vertices of a graph are used to defend it against attacks
on either its vertices or its edges. Various models for this problem have been
proposed. In this survey we describe a number of these models with particular
attention to the case when the attack sequence is infinitely long and the
guards must induce some particular configuration before each attack, such as a
dominating set or a vertex cover. Results from the literature concerning the
number of guards needed to successfully defend a graph in each of these
problems are surveyed.

\end{abstract}

\noindent\textbf{Keywords:\hspace{0.1in}}graph protection, eternal dominating
set, eternal vertex cover

\noindent\textbf{AMS Subject Classification Number 2000:\hspace{0.1in}}05C69

\section{Introduction}

Graph protection involves the placement of mobile guards on the vertices of a
graph to protect its vertices and edges against single or sequences of attacks
and has its historical roots in the time of the ancient Roman Empire. The
modern study of graph protection was initiated in the late twentieth century
by the appearance of four publications in quick succession that referred to
the military strategy of Emperor Constantine (Constantine The Great, 274-337 AD).

The seminal paper is Ian Stewart's \textquotedblleft Defend the Roman
Empire!\textquotedblright\ in \emph{Scientific American}, December 1999
\cite{Stewart}, which contains a reply to C. S. ReVelle's \textquotedblleft
Can you protect the Roman Empire?\textquotedblright, \emph{Johns Hopkins
Magazine}, April 1997 \cite{ReVelle}, and which is based on ReVelle and K. E.
Rosing's \textquotedblleft Defendens Imperium Romanum: A Classical Problem in
Military Strategy\textquotedblright\ in \emph{American Mathematical Monthly},
August -- September 2000~\cite{ReVelleRos}. ReVelle's work \cite{ReVelle} in
turn is a response to the paper \textquotedblleft\ Graphing' an Optimal Grand
Strategy\textquotedblright\ by J. Arquilla and H. Fredricksen \cite{AF}, which
appeared in \emph{Military Operations Research} in 1995 and which is the
oldest reference we could find that places the strategy of Emperor Constantine
in a mathematical setting.

According to ancient history -- some say mythology -- Rome was founded by
Romulus and Remus in 760 -- 750 BC on the banks of the Tiber in central Italy.
It was a country town whose power gradually grew until it was the centre of a
large empire. In the third century AD Rome dominated not only Europe, but also
North Africa and the Near East. The Roman army at that time was strong enough
to use a \emph{forward defense} strategy, deploying an adequate number of
legions to secure on-site every region throughout the empire. However, the
Roman Empire's power was greatly reduced over the following hundred years. By
the fourth century AD only twenty-five legions of the Roman army were
available, which made a forward defense strategy no longer feasible.%
\begin{figure}[ptb]%
\centering
\includegraphics[
height=4.1753in,
width=5.7648in
]%
{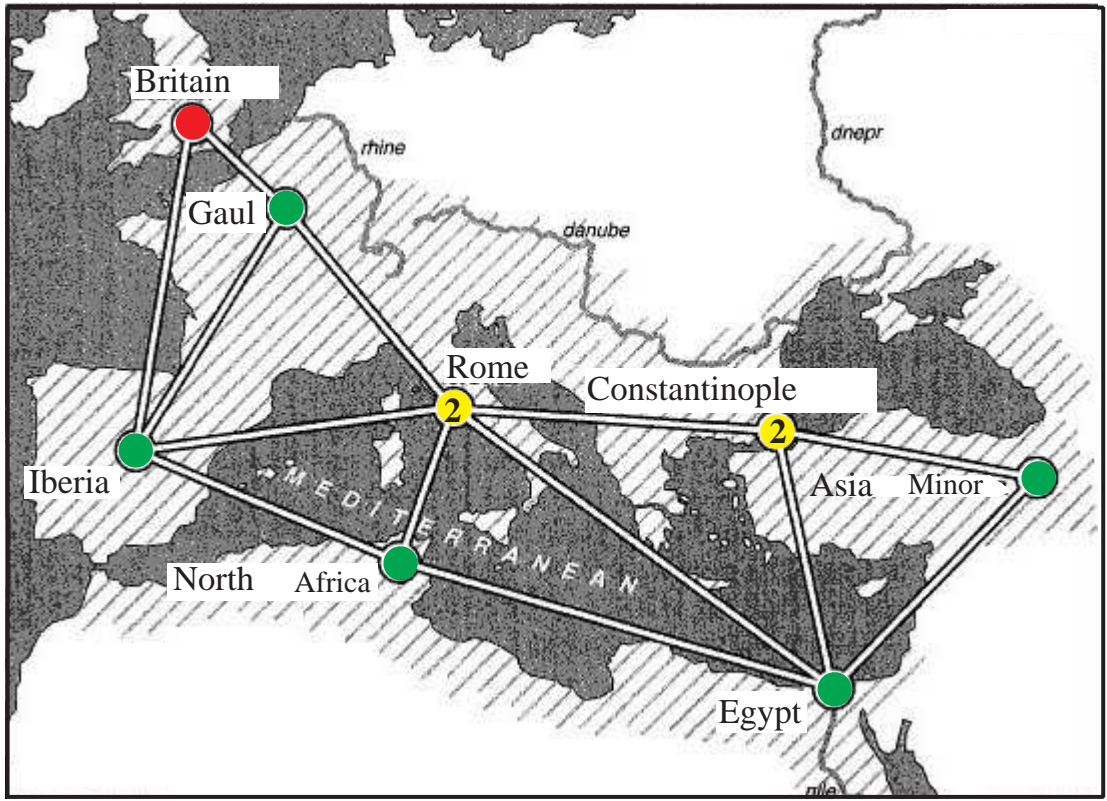}%
\caption{The Roman Empire, fourth century AD}%
\label{Map}%
\end{figure}

According to E. N. Luttwak, \emph{The Grand Strategy of the Roman Empire}, as
cited in \cite{ReVelleRos}, to cope with the reducing power of the Empire,
Constantine devised a new strategy called a \emph{defense in depth} strategy,
which used local troops to disrupt invasion. He deployed mobile Field Armies
(FAs), units of forces consisting of roughly six legions powerful enough to
secure any one of the regions of the Roman Empire, to stop the intruding
enemy, or to suppress insurrection. By the fourth century AD there were only
four FAs available for deployment, whereas there were eight regions to be
defended (Britain, Gaul, Iberia, Rome, North Africa, Constantinople, Egypt and
Asia Minor) in the empire. See Figure \ref{Map}. An FA was considered capable
of deploying to protect an adjacent region only if it moved from a region
where there was at least one other FA to help launch it. The challenge that
Constantine faced was to position four FAs in the eight regions of the empire.
Consider a region to be \emph{secured} if it has one or more FAs stationed in
it already, and \emph{securable} if an FA can reach it in one step.
Constantine decided to place two FAs in Rome and another two FAs in
Constantinople, making all regions either secured or securable -- with the
exception of Britain, which could only be secured after at least four
movements of FAs.

It is mentioned in \cite{AF, ReVelleRos, Stewart} that Constantine's
\textquotedblleft defense in depth\textquotedblright\ strategy was adopted
during World War II by General Douglas MacArthur. When conducting military
operations in the Pacific theatre he pursued a strategy of \textquotedblleft
island-hopping\textquotedblright\ -- moving troops from one island to a nearby
one, but only when he could leave behind a large enough garrison to keep the
first island secure. The efficiency of Constantine's strategy under different
criteria, and ways in which it can be improved, were also discussed in these
three articles.

Constantine's strategy is now known in domination theory as \textbf{Roman
domination}. A \emph{Roman dominating function} on a graph $G=(V,E)$ is a
function $f:V\rightarrow\{0,1,2\}$ satisfying the condition that every vertex
$u$ with $f(u)=0$ is adjacent to at least one vertex $v$ with $f(v)=2$.
\textbf{Weak Roman domination}, an alternative defense strategy that can be
used if defense units can move without another unit being present, was
introduced in \cite{HH2003}. A function $f:V\rightarrow\{0,1,2\}$ is a
\emph{weak Roman dominating function of }$G$ if each vertex $u$ with $f(u)=0$
is adjacent to a vertex $v$ with $f(v)>0$ such that the function $f^{\prime
}=(f-\{(v,f(v)),(u,0)\})\cup\{(v,f(v)-1),(u,1)\}$ also has the property that
each vertex labelled 0 is adjacent to a vertex with positive label.
\textbf{Secure domination} is a defense strategy that can be used when it is
not possible or desirable to station two defense units at the same location. A
\emph{secure dominating function }is a weak Roman dominating function $f$ such
that $\{v\in V:f(v)=2\}=\varnothing$. In this case the set $\{v\in V:f(v)=1\}$
is a \emph{secure dominating set }of~$G$.

A full discussion of Roman domination, weak Roman domination and secure
domination is beyond the scope of this survey. Publications covering these
topics and their variations are given in the bibliography after the list of
references. Here we focus on securing the vertices and edges of graphs against
infinite sequences of attacks, executed one at a time, by stationing defense
units, henceforth called guards, at the vertices of the graph. At most one
guard is stationed at each vertex, and guards that move in response to an
attack do not return to their original positions before facing another attack.
We refer to such models as \emph{eternal}, as they can be thought of as
protecting a graph for eternity. A number of different eternal protection
models have been studied. We introduce them in the next section.

\section{Graph protection models}

\label{Sec_Models}A \emph{dominating set }of a graph $G=(V,E)$ is a set
$D\subseteq V$ such that each vertex in $V-D$ is adjacent to a vertex in $D$.
The minimum cardinality amongst all dominating sets of $G$ is the
\emph{domination number} $\gamma(G)$. By imposing conditions on the subgraph
$G[D]$ of $G$ induced by $D$ we obtain several varieties of dominating sets
and their associated parameters. For example, if $G[D]$ is connected, then $D$
is a \emph{connected dominating set }and the corresponding parameter is the
\emph{connected} \emph{domination number} $\gamma_{c}(G)$, and if $G[D]$ has
no isolated vertices, then $D$ is a \emph{total dominating set }and the
minimum cardinality amongst all total dominating sets is the \emph{total
domination number} $\gamma_{t}(G)$. Only connected graphs have connected
dominating sets, and only graphs without isolated vertices have total
dominating sets. Domination theory can be considered the precursor to the
study of graph protection: one may view a dominating set as an immobile set of
guards protecting a graph. A thorough survey of domination theory can be found
in \cite{HHS}.

A \emph{vertex cover} of $G$ is a set $C\subseteq V$ such that each edge of
$G$ is incident with a vertex in $C$. The minimum cardinality of a vertex
cover of $G$ is the \emph{vertex cover number} (also sometimes called the
\emph{vertex covering number})\emph{ }$\tau(G)$ of $G$. An \emph{independent
set} of $G$ is a set $I\subseteq V$ such that no two vertices in $I$ are
adjacent. The maximum cardinality amongst all independent sets is the
\emph{independence number} $\alpha(G)$. The independence number of $G$ equals
the clique number\emph{ }$\omega(\overline{G})$ of the complement
$\overline{G}$ of $G$. It is well known that $\alpha(G)+\tau(G)=n$ for all
graphs $G$ of order $n$ (see e.g. \cite[p. 241]{CL}). A \emph{matching }in $G$
is a set of edges, no two of which have a common end-vertex. The
\emph{matching number }$m(G)$ is the maximum cardinality of a matching of $G$.
It is also well known that $\tau(G)\geq m(G)$ for all graphs, and that
equality holds for bipartite graphs. The latter result is known as
\textbf{K\"{o}nig's theorem} (see e.g. \cite[Theorem 9.13]{CL}).


Let $\{D_{i}\}$, $D_{i}\subseteq V$, $i\geq1$, be a collection of sets of
vertices of the same cardinality, with one guard located on each vertex of
$D_{i}$. Each protection strategy can be modeled as a two-player game between
a \emph{defender} and an \emph{attacker}: the defender chooses $D_{1}$ as well
as each $D_{i}$, $i>1$, while the attacker chooses the locations of the
attacks $r_{1},r_{2},\ldots$. Each attack is dealt with by the defender by
choosing the next $D_{i}$ subject to some constraints that depend on the
particular game. The defender wins the game if they can successfully defend
any sequence of attacks, subject to the constraints of the game described
below; the attacker wins otherwise.

We say that a vertex (edge) is \emph{protected} if there is a guard on the
vertex or on an adjacent (incident) vertex. A vertex $v$ is \emph{occupied} if
there is a guard on $v$, otherwise $v$ is \emph{unoccupied}. An attack is
\emph{defended} if a guard moves to the attacked vertex (across the attacked edge).

For the \textbf{eternal domination problem}, each $D_{i}$, $i\geq1$, is
required to be a dominating set, $r_{i}\in V$ (assume without loss of
generality $r_{i}\notin D_{i}$), and $D_{i+1}$ is obtained from $D_{i}$ by
moving one guard to $r_{i}$ from an adjacent vertex $v\in D_{i}$. If the
defender can win the game with the sets $\{D_{i}\}$, then each $D_{i}$ is an
\emph{eternal dominating set}. The size of a smallest eternal dominating set
of $G$ is the \emph{eternal domination number} $\gamma^{\infty}(G)$. This
problem was first studied by Burger et al. in \cite{BCG2} and will sometimes
be referred to as the \emph{one-guard moves} model.

For the \textbf{m-eternal dominating set problem}, each $D_{i}$, $i\geq1$, is
required to be a dominating set, $r_{i}\in V$ (assume without loss of
generality $r_{i}\notin D_{i}$), and $D_{i+1}$ is obtained from $D_{i}$ by
moving guards to neighboring vertices. That is, each guard in $D_{i}$ may move
to an adjacent vertex, as long as one guard moves to $r_{i}$. Thus it is
required that $r_{i}\in D_{i+1}$. The size of a smallest \emph{m-eternal
dominating set} (defined similar to an eternal dominating set) of $G$ is the
\emph{m-eternal domination number} $\gamma_{\mathrm{m}}^{\infty}(G)$. This
\textquotedblleft multiple guards move\textquotedblright\ version of the
problem was introduced by Goddard, Hedetniemi and Hedetniemi \cite{GHH}. It is
also called the \textquotedblleft all-guards move\textquotedblright\ model. It
is clear that $\gamma^{\infty}(G)\geq\gamma_{\mathrm{m}}^{\infty}(G)\geq
\gamma(G)$ for all graphs $G$.

As for dominating sets, we obtain variations on the above-mentioned protection
models by imposing conditions on $G[D_{i}]$. Thus we define the \emph{eternal
total }(\emph{connected}, respectively)\emph{ domination number} $\gamma
_{t}^{\infty}(G)$ ($\gamma_{c}^{\infty}(G)$, respectively) and the
\emph{m-eternal total }(\emph{connected}, respectively)\emph{ domination
number} $\gamma_{\mathrm{\operatorname{mt}}}^{\infty}(G)$ ($\gamma
_{\mathrm{mc}}^{\infty}(G)$, respectively) in the obvious way. Eternal total
domination and eternal connected domination were introduced by Klostermeyer
and Mynhardt \cite{KM4}.

For the \textbf{m-eternal vertex covering problem,} each $D_{i}$, $i\geq1$, is
required to be a vertex cover, $r_{i}\in E$, and $D_{i+1}$ is obtained from
$D_{i}$ by moving guards to neighboring vertices; all guards in $D_{i}$ may
move to adjacent vertices provided that one of them moves across edge $r_{i}$
(we assume without loss of generality that one end-vertex of $r_{i}$ is not in
$D_{i}$, otherwise the two guards on the endvertices of $r_{i}$ simply
interchange positions). If the defender can win the game with the sets
$\{D_{i}\}$, then each $D_{i}$ is an \emph{eternal vertex cover}. The size of
a smallest eternal vertex cover of $G$ is the \emph{eternal covering number}
$\tau_{\mathrm{m}}^{\infty}(G)$. The m-eternal vertex covering problem (or
just the eternal vertex covering problem, for simplicity) was introduced by
Klostermeyer and Mynhardt \cite{KM5} and was also studied by Fomin et al. in
\cite{gaspers, gaspers2} and Anderson et al. in \cite{ABC2, ABC}. As in the
case of domination, $\tau_{\mathrm{m}}^{\infty}(G)\geq\tau(G)$ for all graphs
$G$. Also, for any graph $G$ without isolated vertices, $\tau(G)\geq\gamma(G)$
and $\tau_{\mathrm{m}}^{\infty}(G)\geq\gamma_{\mathrm{m}}^{\infty}(G)$.

We discuss these and other related protection models in Sections \ref{Sec_ED}
-- \ref{Sec_other} and present a list of open problems in Section
\ref{Sec_Open}.

We conclude this section with some remarks about the nature of the attack
sequence $\{r_{i}\}$. There are three main ways for the attacker to choose and
reveal $\{r_{i}\}$. Following the notation used for the $k$-server problem
(see Section \ref{Sec_k-server}), they are as follows.

\begin{enumerate}
\item \textbf{{Offline problem}}: the entire sequence $r_{1},r_{2}%
,\ldots,r_{m}$ of attacks is chosen and revealed in advance.

\item \textbf{{Adaptive online problem}}: the sequence of attacks is chosen
and revealed one by one by the attacker alternating with the guard movements
by the defendant. The attacker is called an \emph{adaptive adversary}.

\item \textbf{{Oblivious online problem}}: the sequence of attacks is
constructed in advance by an adversary, but revealed one by one in response to
each guard movement. The adversary in this case is called an \emph{oblivious
adversary}.
\end{enumerate}

The offline problem, even if the finite sequence $r_{1},r_{2},\ldots,r_{m}$ is
repeated indefinitely, is not the same as eternal domination problem. The
minimum number of guards required to defend such a predefined attack sequence
could be strictly less than the eternal domination number. We only consider
this type of attack sequence for the $k$-server problem in Section
\ref{Sec_k-server}. The adaptive online problem is precisely the eternal
domination problem as described above: the location of each attack is chosen
by the attacker depending on the location of the guards at that time. At first
glance, the oblivious online problem appears to be somewhat different from the
adaptive online problem, and to be the same as the original eternal domination
problem described in \cite{BCG2}. However, the defender is required to defend
against \emph{any }attack sequence and has no advance knowledge of the
sequence. Furthermore, one can assume the attacker is aware of the defense
strategy; and so the attacker can predict the moves of the defender unless the
defender employs a randomized strategy. Because randomized strategies are not
relevant for the types of results described in this paper, for our purposes,
the two types of attack models are equivalent. Certainly, the associated
parameters are equal. Randomized strategies are relevant when one asks
questions that might concern the number of (expected) moves before some
configuration is reached, for example.

\section{Definitions}

The \emph{open} and \emph{closed neighbourhoods} of $X\subseteq V$ are
$N(X)=\{v\in V:v$ is adjacent to a vertex in $X\}$ and $N[X]=N(X)\cup X$,
respectively, and $N(\{v\})$ and $N[\{v\}]$ are abbreviated, as usual, to
$N(v)$ and $N[v]$. For any $v\in X$, the \emph{private neighbourhood}
$\operatorname{pn}(v,X)$\emph{ of }$v$\emph{ with respect to }$X$ is the set
of all vertices in $N[v]$ that are not contained in the closed neighbourhood
of any other vertex in $X$, i.e., $\operatorname{pn}(v,X)=N[v]-N[X-\{v\}]$.
The elements of $\operatorname{pn}(v,X)$ are the \emph{private neighbours of
}$v$ \emph{relative to }$X$. The \emph{external private neighbourhood of }$v$
\emph{with respect to }$X$ is the set $\operatorname{epn}%
(v,X)=\operatorname{pn}(v,X)-\{v\}=N(v)-N[X-\{v\}]$.

The \emph{clique covering number }$\theta(G)$ is the minimum number $k$ of
sets in a partition $V=V_{1}\cup\cdots\cup V_{k}$ of $V$ such that each
$G[V_{i}]$ is complete. Hence $\theta(G)$ equals the chromatic number
$\chi(\overline{G})$ of the complement $\overline{G}$ of $G$. Since
$\chi(G)=\omega(G)$ if $G$ is perfect, and $G$ is perfect if and only if
$\overline{G}$ is perfect \cite[p. 203]{CL}, $\alpha(G)=\theta(G)$ for all
perfect graphs.

The \emph{circulant graph} $C_{n}[a_{1},...,a_{k}]$, where $1\leq a_{1}%
\leq\cdots\leq a_{k}\leq\left\lfloor \frac{n}{2}\right\rfloor $, is the graph
with vertex set $\{v_{0},...,v_{n-1}\}$ such that $v_{i}$ and $v_{j}\ $are
adjacent if and only if $i-j\equiv\pm a_{\ell}\ (\operatorname{mod}\ n)$ for
some $\ell\in\{1,...,k\}$.

The Cartesian product of two graphs $G$ and $H$ is denoted $G\boksie H$; a
definition can be found in \cite{HHS}.

\section{Eternal domination}

\label{Sec_ED} The eternal domination problem was first studied by Burger et
al. \cite{BCG2} in 2004 where it was called infinite order domination. That
paper, and this section, consider the one-guard moves model. Shortly
thereafter, Goddard et al. published a second paper on the subject where they
called it \emph{eternal security} \cite{GHH}.

Consider an eternal dominating set $D$ of a graph $G$. A necessary condition
for a guard on $D$ to defend a neighbouring vertex in a winning strategy is
given below.

\begin{proposition}
\label{Prop_winning}Let $D$ be an eternal dominating set of a graph $G$. If a
guard on $v\in D$ can move to a vertex $u\in V-D$ in a winning strategy, then
$\operatorname{pn}(v,D)\cup\{u\}$ induces a clique.
\end{proposition}

\noindent\textbf{Proof}.\hspace{0.1in}Suppose the guard $g$ on $v$ moves to
$u$ in a winning strategy. If the next attack is at $x\in\operatorname{pn}%
(v,D)$, $g$ moves to $x$, as it is the only guard that protects $x$. Since
this holds whether $u\in\operatorname{pn}(v,D)$ or not, $\operatorname{pn}%
(v,D)\cup\{u\}$ induces a clique.\hfill$\blacksquare$

\bigskip

\begin{figure}[ptb]%
\centering
\includegraphics[
height=1.663in,
width=4.4521in
]%
{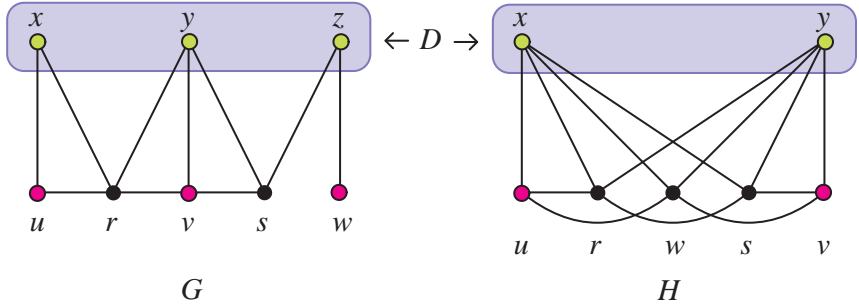}%
\caption{In $G$, $y$ does not defend $r$, and $D$ is not an eternal dominating
set of $H$}%
\label{Fig_Protect1}%
\end{figure}

The converse of Proposition \ref{Prop_winning} is not true. Consider the graph
$G$ in Figure \ref{Fig_Protect1}. The set $D=\{x,y,z\}$ is an eternal
dominating set of $G$ in which the guard on $x$ ($y$, $z$) defends $\{x,u,r\}$
($\{y,v,s\}$, $\{z,w\}$). Also, $\operatorname{pn}(y,D)=\{y,v\}$ and
$G[\{y,v,r\}]$ is a clique. Suppose, however, the guard on $y$ moves to $r$.
If the next attack is at $s$, then only $z$ has a guard adjacent to $s$. But
moving this guard to $s$ leaves $w$ unprotected. In the graph $H$ in
Figure~\ref{Fig_Protect1}, $D=\{x,y\}$ is not an eternal dominating set, even
though $\operatorname{pn}(x,D)\cup\{r\}$, $\operatorname{pn}(x,D)\cup\{w\}$,
$\operatorname{pn}(y,D)\cup\{w\}$, $\operatorname{pn}(y,D)\cup\{s\}$ all
induce cliques: first attack $r$; without loss of generality, the guard on $x$
moves there. Now attack $s$. If the guard on $y$ moves there, then $w$ is not
protected; if the guard on $r$ moves there, then $u$ is not protected.%

\subsection{Bounds for the eternal domination number\label{Sec_Bounds}}

As first observed by Burger et al. \cite{BCG2}, it does not take much
imagination to see that $\gamma^{\infty}$ lies between the independence and
clique covering numbers.

\begin{Fact}
\label{FactED_Bound}For any graph $G$, $\alpha(G)\leq\gamma^{\infty}%
(G)\leq\theta(G)$.
\end{Fact}

\noindent\textbf{Proof}.\hspace{0.1in}To see the lower bound, consider a
sequence of consecutive attacks at the vertices of a maximum independent set.
To see the upper bound, observe that a single guard can defend all vertices of
a clique. \hfill$\blacksquare$

\bigskip

Since $\alpha(G)=\theta(G)$ for perfect graphs, the bounds in Fact
\ref{FactED_Bound} are tight for perfect graphs. A topic that has received
much attention is finding classes of non-perfect graphs that satisfy one of
the bounds in Fact \ref{FactED_Bound}. Before proceeding, we should point out
that the independence number, eternal domination number, and clique-covering
number can vary widely.

\begin{theorem}
\label{ThmDifference}\emph{\cite{KM2}}\hspace{0.1in}For any positive integers
$c$ and $d$ there exists a connected graph $G$ such that $\alpha
(G)+c\leq\gamma^{\infty}(G)$ and $\gamma^{\infty}(G)+d\leq\theta(G)$.
\end{theorem}

Let $C_{n}^{k}$ denote the $k^{th}$ power (see \cite[p. 105]{CL}) of the cycle
of order $n$, where $2k+1<n$.

\begin{theorem}
If $G$ is a graph in one of the following classes, then $\gamma^{\infty
}(G)=\theta(G)$.

\begin{enumerate}
\item[\emph{(a)}] \emph{\cite{BCG2}}\hspace{0.1in}Perfect graphs.

\item[\emph{(b)}] \emph{\cite{BCG2}}\hspace{0.1in}Any graph $G$ such that
$\theta(G)\leq3$.

\item[\emph{(c)}] \emph{\cite{KM2}}$\hspace{0.1in}C_{n}^{k}$ and
$\overline{C_{n}^{k}}$, for all $k\geq1$, $n\geq3$.

\item[\emph{(d)}] \emph{\cite{regan}}\hspace{0.1in}Circular-arc graphs
(intersection graphs of a family of arcs of a circle).

\item[\emph{(e)}] \emph{\cite{ABB}}$\hspace{0.1in}K_{4}$-minor-free graphs
(a.k.a. series-parallel graphs, see e.g. \cite[p. 336]{West} for definition).

\item[\emph{(f)}] \emph{\cite{ABB}}$\hspace{0.1in}C_{m}\boksie C_{n}$;
$P_{m}\boksie C_{n}$.
\end{enumerate}
\end{theorem}

Goddard et al. \cite{GHH} showed that
\begin{equation}
\text{if }\alpha(G)=2,\text{ then }\gamma^{\infty}(G)\leq3.
\label{eq_GHH_Bound}%
\end{equation}
The Mycielski construction (see \cite[p. 203]{CL}) yields triangle-free
$k$-chromatic graphs for arbitrary $k$. The complements of these graphs have
$\alpha=2$ and $\theta=k$, and hence are examples of graphs with small eternal
domination numbers and large clique covering numbers. The Gr\"{o}tzsch graph
is the smallest $4$-chromatic triangle-free graph, and its complement is the
smallest known graph with $\gamma^{\infty}<\theta$. Goddard et al. \cite{GHH}
also gave the first example of a graph $G$ with $\alpha(G)<\gamma^{\infty
}(G)<\theta(G)$: the circulant graph $C_{18}[1,3,8]$, which satisfies
$\alpha=6$, $\gamma^{\infty}=8$ and $\theta=9$.

Klostermeyer and MacGillivray \cite{klos1} proved the existence of graphs with
$\gamma^{\infty}=\alpha$ and whose clique covering number is either equal to
two (if $\alpha=2$) or arbitrary otherwise. Their proof rests (i.~a.) on the
observation that if $H$ is an induced subgraph of $G$ and $\pi$ is any of the
parameters $\alpha$, $\gamma^{\infty}$, $\theta$, then $\pi(H)\leq\pi(G)$.
This is trivially true for $\alpha$ and $\theta$. To see that it is true for
$\gamma^{\infty}$, note that a sequence of attacks on $G$ but restricted to
$H$ requires $\gamma^{\infty}(H)$ guards, hence $\gamma^{\infty}(G)\geq
\gamma^{\infty}(H)$.

\begin{theorem}
\label{Thm_equal_LB}\emph{\cite{klos1}}

\begin{enumerate}
\item[\emph{(a)}] If $\alpha(G)=\gamma^{\infty}(G)=2$, then $\theta(G)=2$.

\item[\emph{(b)}] For all integers $k\geq a\geq3$ there exists a connected
graph $G$ such that $\alpha(G)=\gamma^{\infty}(G)=a$ and $\theta(G)=k$.
\end{enumerate}
\end{theorem}

\noindent\textbf{Proof}.\hspace{0.1in}(a)\hspace{0.1in}The statement is
clearly true for graphs of order three and four. Assume it to be true for all
graphs of order less than $n$, where $n\geq5$, and let $G$ be an $n$-vertex
graph with $\alpha(G)=\gamma^{\infty}(G)=2$. Let $u$ and $v$ be nonadjacent
vertices of $G$. After consecutive attacks on $u$ and $v$, both these vertices
are occupied. By Proposition \ref{Prop_winning}, $\operatorname{pn}%
(v,\{u,v\})=V-N[u]$ and $\operatorname{pn}(u,\{u,v\})=V-N[v]$ induce cliques.
Let $W$ and $Y$ be the sets of all vertices in $N(u)\cap N(v)$ that are
defended by the guard on $v$ and the guard on $u$, respectively. By
Proposition \ref{Prop_winning}, each $w\in W$ ($y\in Y$, respectively) is
adjacent to each vertex in $V-N[u]$ ($V-N[v]$, respectively).

Let $H^{\prime}\cong G[N(u)\cap N(v)]$ and $H\cong G[V(H^{\prime}%
)\cup\{u,v\}]$. (Possibly $H=G$.) If $H^{\prime}$ is complete, then
$(V-N[u])\cup W$ and $(V-N[u])\cup Y$ induce cliques that contain $V(G)$.
Hence suppose $H^{\prime}$ is not complete. Then $\alpha(H^{\prime})\geq2$ and
so $\gamma^{\infty}(H^{\prime})\geq2$. Since $H^{\prime}$ is an induced
subgraph $H$, which is an induced subgraph of $G$, it follows that
$\alpha(H)=\gamma^{\infty}(H)=2$ and $\alpha(H^{\prime})=\gamma^{\infty
}(H^{\prime})=2$. By the induction hypothesis, $\theta(H^{\prime})=2$ and so
$\theta(H)=2$. Partition $V(H)$ into the cliques $C_{u},C_{v}$, where $u\in
V(C_{u})$, $v\in V(C_{v})$. Clearly, $V(C_{u})-\{u\}\subseteq W$ and
$V(C_{v})-\{v\}\subseteq Y$. Therefore $(V-N[v])\cup V(C_{u})$ and
$(V-N[u])\cup V(C_{v})$ induce a clique partition of $G$.

\bigskip

\noindent(b)\hspace{0.1in}Let $H$ be the complement of a triangle-free
$k$-chromatic graph, $k\geq3$. Then $\alpha(H)=2$ and $\theta(H)=k$. By (a),
$\gamma^{\infty}(H)\geq3$, and thus by (\ref{eq_GHH_Bound}), $\gamma^{\infty
}(H)=3$. Add $a$ new vertices $v_{1},...,v_{a}$, joining each $v_{i}$ to each
vertex of $H$ to form the graph $G$. Then $\alpha(G)=a$, and, since $a\leq k$,
$\theta(G)=k$. Place a guard on each $v_{i}$, $i>3$; these guards never move.
The remaining three guards protect $H$ and $v_{1},v_{2},v_{3}$ according to
the strategy for $H$; when $v_{i}$ is attacked, any guard moves there, and
returns to $H$ when required.\hfill$\blacksquare$

\bigskip

Goddard et al. \cite{GHH} asked whether the eternal domination number can be
bounded by a constant times the independence number. That this is impossible
in general follows from the next two theorems. One of the main results on
eternal domination is the following upper bound, due to Klostermeyer and
MacGillivray \cite{KM}.

\begin{theorem}
\label{ThmKM}\emph{\cite{KM}\hspace{0.1in}}For any graph $G$,%
\[
\gamma^{\infty}(G)\leq\binom{\alpha(G)+1}{2}.
\]

\end{theorem}

\noindent\textbf{Proof}.\hspace{0.1in}Assume $|V|>{\binom{\alpha+1}{2}}$,
otherwise we are done. Consider pairwise disjoint sets $S_{\alpha}%
,S_{\alpha-1},\ldots,S_{1}$, where $S_{\alpha}$ is a maximum independent set
of $G$ and, for $i=\alpha-1,\alpha-2,\ldots,1$, the set $S_{i}$ is either
empty or an independent set of size $i$. Other than $S_{\alpha}$, no $S_{i}$
needs to be a maximal independent set. Among all collections of such sets, we
choose one such that $|\bigcup_{i=1}^{\alpha}S_{i}|$ is maximum. Since
$|V|>{\binom{\alpha+1}{2}}$, the set $S_{1}\neq\varnothing$. Let
$D=\bigcup_{i=1}^{\alpha}S_{i}$ and note that $|D|\leq\binom{\alpha(G)+1}{2}$.
We describe a defense strategy $\bigstar$ which shows that $D$ is an eternal
dominating set of $G$.

\begin{itemize}
\item[$\bigstar$] Whenever there is an attack at a vertex $v\notin D$, a guard
on a vertex $u$ from the set $S_{t}$ with the smallest subscript among those
with a vertex adjacent to $v$ moves to $v$. Such a set $S_{t}$ exists because
$S_{\alpha}$ is a dominating set (as it is a maximum independent set).
\end{itemize}

The key technical part of the proof is to show that $(D-\{u\})\cup\{v\}$ can
be partitioned into disjoint independent sets with the same properties as the
sets $\{S_{i}\}$. There are two cases.

If $S_{t}^{\prime}=(S_{t}-\{u\})\cup\{v\}$ is an independent set, then
replacing $S_{t}$ by $S_{t}^{\prime}$ yields another collection of disjoint
independent sets as desired. Otherwise, $v$ is adjacent to at least two
vertices in $S_{t}$ and $t>1$. Let $r$ be the greatest integer less than $t$
such that $S_{r}\neq\varnothing$. We show that $r=t-1$.

Suppose $r\leq t-2$. Then $S_{r+1}=\varnothing$. By definition of $t$, no
vertex in $S_{r}$ is adjacent to $v$, hence $S_{r}\cup\{v\}$ is an independent
set of cardinality $r+1$. The collection of independent sets obtained by
replacing $S_{r+1}$ by $S_{r}\cup\{v\}$ and $S_{r}$ by $\varnothing$
contradicts the maximality of $|\bigcup_{k=1}^{\alpha}S_{k}|$. Hence $r=t-1$.

Replacing $S_{t}$ by $S_{t-1}\cup\{v\}$ and $S_{t-1}$ by $S_{t}-\{u\}$ gives
another collection of independent sets with the desired properties. Thus we
may repeat $\bigstar$ indefinitely to protect $G$ against any sequence of
attacks.\hfill$\blacksquare$

\bigskip

Goldwasser and Klostermeyer \cite{GK} showed that this bound is sharp for
certain graphs. Specifically, let $G(n,k)$ be the graph with vertex set
consisting of the set of all $k$-subsets of an $n$-set and where two vertices
are adjacent if and only if their intersection is nonempty (thus $G(n,k)$ is
the complement of a Kneser graph).

\begin{theorem}
\label{ThmGK}\emph{\cite{GK}\hspace{0.1in}}For each positive integer $t$, if
$k$ is sufficiently large, then the graph $G(kt+k-1,k)$ has eternal domination
number ${\binom{t+1}{2}}$.
\end{theorem}

Regan \cite{regan} found another graph for which the bound is sharp: the
circulant graph $C_{22}[1,2,4,5,9,11]$. Theorems \ref{ThmKM} and \ref{ThmGK}
show that it is impossible to find a constant $c$ such that $\gamma^{\infty
}(G)\leq c\alpha(G)$ for all graphs $G$. It would be of interest to find other
graphs where the bound is sharp.

As shown by Klostermeyer and MacGillivray \cite{KM}, the graph $G$ obtained by
joining a new vertex to $m$ disjoint copies of $C_{5}$ satisfies
$\alpha(G)=2m$ and $\gamma^{\infty}(G)=3m$, that is, $\gamma^{\infty
}(G)/\alpha(G)=\frac{3}{2}$. This result and Theorem \ref{Thm_equal_LB} can be
placed in a more general setting, as explained in \cite{KM2}.

A triple $(a,g,t)$ of positive integers is called \emph{realizable }if there
exists a connected graph $G$ with $\alpha(G)=a$, $\gamma^{\infty}(G)=g$ and
$\theta(G)=t$. Theorem \ref{ThmKM} shows that no triple with $g>\binom{a+1}%
{2}$ is realizable. The following theorem, stated in \cite{KM2}, provides a
partial solution to the question of which triples are realizable.

\begin{theorem}
\label{ThmTriples}Let $(a,g,t)$ be a triple of positive integers such that
$a\leq g\leq t$.

\begin{enumerate}
\item[\emph{(a)}] The only realizable triple with $a=1$ is $(1,1,1)$.

\item[\emph{(b)}] \emph{\cite{BCG2, GHH, klos1}}\hspace{0.1in}The only
realizable triples with $a=2$ are $(2,2,2)$ and $(2,3,t)$, $t\geq3$.

\item[\emph{(c)}] \emph{\cite{BCG2, klos1, KM2}\hspace{0.1in}}For all integers
$a$, $g$ and $t$ with $3\leq a\leq g\leq\frac{3}{2}a$ and $g\leq t$, the
triple $(a,g,t$) is realizable.
\end{enumerate}
\end{theorem}

The circulant $C_{21}[1,3,8]$, which satisfies $\gamma^{\infty}/\alpha
=\frac{10}{6}$ (see \cite{GHH}), shows that Theorem \ref{ThmTriples} does not
characterize realizable triples.

\subsection{The $k$-server problem}

\label{Sec_k-server}We briefly mention the $k$-server problem, which is
related to the eternal domination problem. The $k$-server problem is an
algorithmic problem set in the more general framework of metric spaces, but
often focused on graphs. It was defined in \cite{mms} as follows. There are
$k$ mobile servers (or guards) located at vertices of a graph. In response to
an attack on an unoccupied vertex $r_{i}$, a server must move to $r_{i}$. The
objective is to minimize the total distance travelled by all the servers over
the sequence of attacks. The three main variations of the problem are (1) the
offline problem, (2) the adaptive online problem and (3) the oblivious online
problem, as described in Section \ref{Sec_Models}.

A simple polynomial time algorithm using dynamic programming can compute the
optimal solution for the offline problem \cite{mms}. A faster algorithm is
given in \cite{ckpv}.

Koutsoupias and Papadimitriou \cite{kout} proved that a simple algorithm known
as the work-function algorithm is $2k-1$ competitive. In other words, the
distance the servers travel using the work function algorithm is at most
$2k-1$ times the distance they would travel using any other algorithm,
including an optimal algorithm that knew the entire attack sequence in
advance, over all attack sequences. It is a famous conjecture in computer
science that the work function algorithm is $k$-competitive and that this
would be best possible.

A key difference between problems (2) and (3) is that a randomized algorithm
can be effective in problem (3). Since an oblivious adversary cannot adapt the
attack sequence to the moves of the algorithm, by using randomization an
algorithm may be able to effectively prevent an adversary from constructing a
costly attack sequence. A famous result from \cite{ms} is an $H_{k}%
$-competitive algorithm for the problem of $k$ servers on a complete graph
with $k+1$ vertices, where $H_{k}$ is the $k^{th}$ harmonic number. This
result is known to be optimal.

\section{m-Eternal domination}

\label{Sec_mED}

As mentioned in Section \ref{Sec_Models}, m-eternal dominating sets are
defined similar to eternal dominating sets, except that when an attack occurs,
each guard is allowed to move to a neighbouring vertex to either defend the
attacked vertex or to better position themselves for the future. This model
was introduced by Goddard et al. \cite{GHH}. As stated above, we refer to this
as the \textquotedblleft all-guards move\textquotedblright\ model of eternal domination.

Goddard et al. \cite{GHH} determine $\gamma_{\mathrm{m}}^{\infty}(G)$ exactly
for complete graphs, paths, cycles, and complete bipartite graphs. They also
obtained the following fundamental bound.

\begin{theorem}
\label{GHH1} \emph{\cite{GHH}}\hspace{0.1in}For all graphs $G$, $\gamma
(G)\leq\gamma_{\mathrm{m}}^{\infty}(G)\leq\alpha(G)$.
\end{theorem}

\noindent\textbf{Outline of proof}.\hspace{0.1in}The left inequality is
obvious. The right inequality is proved by induction on the order of $G$, the
result being easy to see for small graphs. If $G$ has a vertex $v$ that is not
contained in any maximum independent set, then $v$ is adjacent to at least two
vertices of each maximum independent set of $G$. Therefore $\alpha
(G-N[v])\leq\alpha(G)-2$. Hence (by induction) $G-N[v]$ can be protected by
$\alpha(G)-2$ guards. Since $K_{1,\deg(v)}$ is a spanning subgraph of
$G[N[v]]$, $G[N[v]]$ can be protected by two guards. It follows that
$\gamma_{\mathrm{m}}^{\infty}(G)\leq\alpha(G)$.

If each vertex of $G$ is contained in a maximum independent set, place a guard
on each vertex of a maximum independent set $M$. Defend an attack on $v\in
V(G)-M$ by moving all guards to a maximum independent set $M_{v}$ containing
$v$. This is possible since Hall's Marriage Theorem ensures that there is a
matching between $M_{v}$ and $M$. \hfill$\blacksquare$

\bigskip

Theorem \ref{GHH1} places $\gamma_{\mathrm{m}}^{\infty}$ nicely in the chain
\[
\gamma(G)\leq\gamma_{\mathrm{m}}^{\infty}(G)\leq\alpha(G)\leq\gamma^{\infty
}(G)\leq\theta(G).
\]

Goddard et al.~also claim that $\gamma_{\mathrm{m}}^{\infty}(G)=\gamma(G)$ for
all Cayley graphs~$G$. This claim, however, is false, as is shown in the
recent paper \cite{BSL} by Graga, de Souza and Lee. By computing $\gamma(G)$
and $\gamma_{\mathrm{m}}^{\infty}(G)$ for 7871 Cayley graphs of non-abelian
groups, they found 61 connected Cayley graphs $G$ such that $\gamma
_{\mathrm{m}}^{\infty}(G)=\gamma(G)+1$. For all other connected Cayley graphs
they investigated, $\gamma_{\mathrm{m}}^{\infty}(G)=\gamma(G)$.

The upper bound in Theorem \ref{GHH1} is not tight in general. For example,
$K_{1,m}$ has independence number $m$ and can be defended with just two guards
in this model. But equality holds for many graphs, such as $K_{n}$, $C_{n}$,
and $P_{2}\boksie P_{3}$, just to name a few. By a careful analysis of the
clique structure, it was shown in \cite{Braga2} that $\gamma_{\mathrm{m}%
}^{\infty}(G)=\alpha(G)$ for all proper-interval graphs (a subclass of perfect
graphs). Characterizing graphs with m-eternal domination number equal to the
bounds in Theorem \ref{GHH1} remains open, as mentioned in Section
\ref{SecOpen_m-eternal}. However, trees for which equality holds in the upper
bound, $\alpha$, are characterized by Klostermeyer and MacGillivray
\cite{KM-trees}.

Define a \textit{neo-colonization} to be a partition $\{V_{1},V_{2}%
,\ldots,V_{t}\}$ of $G$ such that each $V_{i}$ induces a connected graph. A
part $V_{i}$ is assigned weight one if it induces a clique, and $1+\gamma
_{c}(G[V_{i}])$ otherwise, where $\gamma_{c}(G[V_{i}])$ is the connected
domination number of the subgraph induced by $V_{i}$. Then $\theta_{c}(G)$ is
the minimum weight of any neo-colonization of $G$.

Goddard et al. \cite{GHH} proved that $\gamma_{\mathrm{m}}^{\infty}%
(G)\leq\theta_{c}(G)\leq\gamma_{c}(G)+1$. Klostermeyer and MacGillivray
\cite{KM2} proved that equality holds in the first inequality for trees.

\begin{theorem}
\emph{\cite{KM2}}\hspace{0.1in}If $T$ is a tree, then $\gamma_{\mathrm{m}%
}^{\infty}(T)=\theta_{c}(T)$.
\end{theorem}

A different upper bound is given in \cite{prince2}. A proof is given below. A
\emph{branch vertex }of a tree\emph{ }is a vertex of degree at least three.

\begin{theorem}
\label{princeupper}\hspace{0.1in}If $G$ is a connected graph of order $n$,
then $\gamma_{\mathrm{m}}^{\infty}(G)\leq\lceil\frac{n}{2}\rceil$.
\end{theorem}

\noindent\textbf{Proof}.\hspace{0.1in}The proof is by induction on $n$, the
result being easy to see for paths and cycles. Let $T$ be a spanning tree of
$G$ with $r\geq1$ branch vertices.

If $T$ has no vertex of degree two, then the subgraph of $T$ induced by the
branch vertices is connected and, by \cite[Theorem 3.7]{CL}, $T$ has at least
$r+2$ leaves. Hence $n\geq2r+2$. Place a guard on each branch vertex and on
one leaf. Whenever an unoccupied leaf $u$ is attacked, guards move so that $u$
and all branch vertices have guards. Hence $\gamma_{\mathrm{m}}^{\infty
}(T)\leq r+1\leq\lceil\frac{n}{2}\rceil$.

If $T$ has a vertex $v$ of degree two, and $N(v)=\{u_{1},u_{2}\}$, then at
least one of the graphs $T-\{vu_{i}\}$, $i=1,2$, has a component of even
order. Let $T_{1}$ be this component and let $T_{2}$ be the other component.
Say $T_{i}$ has order $n_{i}$. By the induction hypothesis, $\gamma
_{\mathrm{m}}^{\infty}(T_{1})\leq\frac{n_{1}}{2}$ and $\gamma_{\mathrm{m}%
}^{\infty}(T_{2})\leq\lceil\frac{n_{2}}{2}\rceil$. It follows that
$\gamma_{\mathrm{m}}^{\infty}(T)\leq\lceil\frac{n}{2}\rceil$ and therefore
$\gamma_{\mathrm{m}}^{\infty}(G)\leq\gamma_{\mathrm{m}}^{\infty}(T)\leq
\lceil\frac{n}{2}\rceil$.\hfill$\blacksquare$

\bigskip

The bound in Theorem \ref{princeupper} is exact for the coronas of all graphs
because they have domination numbers equal to half their order.

It is not hard to see that for many all-guards move models, the associated
parameter is bounded above by $2\gamma$.

\begin{proposition}
\label{twice} For any connected graph $G$, $\gamma_{\mathrm{m}}^{\infty
}(G)\leq2\gamma(G)$, and the bound is sharp for all values of $\gamma(G)$.
\end{proposition}

\noindent\textbf{Proof}.\hspace{0.1in}The result is trivial for $K_{1}$, so
assume $|V(G)|\geq2$. As shown in \cite{BC}, every graph without isolated
vertices has a minimum dominating set in which each vertex has an external
private neighbour. Let $D$ be such a minimum dominating set of $G$. For each
$v\in D$, place a guard at $v$ and at a private neighbour of $v$. This
configuration is an m-eternal dominating set.

To see that the bound is sharp for $\gamma=1$, consider any star with at least
three vertices. For $\gamma=2$, consider $C_{6}$ and let $u$ and $v$ be two
vertices at distance three apart. Add two new internally disjoint $u-v$ paths
of length three to form the graph $G$. Obviously, $\{u,v\}$ is a $\gamma$-set
of $G$. Let $D$ be any dominating set of $G$ with $|D|=3$. Suppose $u\notin
D$. Since $N(u)$ is independent with $|N(u)|=4$, and no two vertices in $N(u)$
have a common neighbour other than $u$, $D$ does not dominate $N(u)$, a
contradiction. Thus $u\in D$ and similarly $v\in D$. Without loss of
generality say $D=\{u,v,w\}$, where $w\in N(u)$. Then $D$ cannot repel an
attack at a vertex in $N(v)-N(w)$. It follows that $\gamma_{\mathrm{m}%
}^{\infty}(G)=4=2\gamma(G)$.

For $\gamma=k\geq3$, consider the cycle $C_{3k}=u_{0},u_{1},...,u_{3k-1}%
,u_{0}$ and the $\gamma$-set $\{u_{0},u_{3},...,u_{3k-3}\}$ of $C_{3k}$. For
each $i=0,...,k-1$, add a new $u_{3i}-u_{3(i+1)(\operatorname{mod}\ 3k)}$ path
of length three to form $G$. Then $\gamma(G)=k$, but it can be shown similar
to the previous case that no set of $2k-1$ vertices eternally protects the
vertices of $G$. \hfill$\blacksquare$

\bigskip

Klostermeyer and MacGillivray \cite{KM-trees} characterized trees for which
equality holds in the following bounds: $\gamma_{\mathrm{m}}^{\infty}%
(T)\leq\gamma_{c}(T)+1$, $\gamma(T)\leq\gamma_{\mathrm{m}}^{\infty}(T)$,
$\gamma_{\mathrm{m}}^{\infty}(T)\leq2\gamma(T)$, and $\gamma_{\mathrm{m}%
}^{\infty}(T)\leq\alpha(T)$.

Grid graphs, i.e. $P_{n}\boksie P_{m}$, are a well-studied class of graphs in
domination theory; see \cite{HHS}. We sometimes refer to $P_{n}\boksie P_{m}$
as the $n\times m$ grid graph. As shown in \cite{GK2}\emph{,} $\gamma
_{\mathrm{m}}^{\infty}(P_{2}\boksie P_{n})=\lceil\frac{2n}{3}\rceil$ for any
$n\geq2$, while $\gamma_{\mathrm{m}}^{\infty}(P_{3}\boksie P_{n})=n$ for
$2\leq n\leq8$. Based on these results, the next two theorems may seem surprising.

\begin{theorem}
\emph{\cite{GK2}} \label{threeby} For $n\geq9$, $\gamma_{\mathrm{m}}^{\infty
}(P_{3}\boksie P_{n})\leq\lceil\frac{8n}{9}\rceil$.
\end{theorem}

\begin{theorem}
\label{3-grid}\emph{\cite{FMV}} For $n>11$, $1+\left\lceil \frac{4n}%
{5}\right\rceil \leq\gamma_{\mathrm{m}}^{\infty}(P_{3}\boksie P_{n}%
)\leq2+\left\lceil \frac{4n}{5}\right\rceil $.
\end{theorem}

Theorem \ref{3-grid} shows that the bound in Theorem \ref{threeby} is not
sharp in general, although it is sharp for $n=9, 10$ for example. It is
conjectured in \cite{GK2} that the lower bound in Theorem\emph{ }\ref{3-grid}
gives the exact value of $\gamma_{\mathrm{m}}^{\infty}(P_{3}\boksie P_{n})$
for $n\geq10$.

Beaton, Finbow and MacDonald \cite{BFM, BFM2} continued the study of m-eternal
domination in grid graphs and obtained the following results.

\begin{theorem}
\emph{\cite{BFM, BFM2}}\ 

\begin{enumerate}
\item[\emph{(a)}] For any $n\in\mathbb{Z}^{+}$, $\gamma_{\mathrm{m}}^{\infty
}(P_{4}\boksie P_{n})=2\left\lceil \frac{n+1}{2}\right\rceil $, with the
exceptions $\gamma_{\mathrm{m}}^{\infty}(P_{4}\boksie P_{2})=3$ and
$\gamma_{\mathrm{m}}^{\infty}(P_{4}\boksie P_{6})=7$.

\item[\emph{(b)}] For any $n\in\mathbb{Z}^{+}$, $\left\lfloor \frac
{10(n+1)}{7}\right\rfloor \leq\gamma_{\mathrm{m}}^{\infty}(P_{6}\boksie
P_{n})\leq\left\lceil \frac{8n}{5}\right\rceil +8$.

\item[\emph{(c)}] $\gamma_{\mathrm{m}}^{\infty}(P_{5}\boksie P_{5})=7$,
$\gamma_{\mathrm{m}}^{\infty}(P_{6}\boksie P_{6})=10$, and $13\leq
\gamma_{\mathrm{m}}^{\infty}(P_{7}\boksie P_{7})\leq14$.
\end{enumerate}
\end{theorem}

Van Bommel and Van Bommel further the study of $5\times n$ grids \cite{vb},
including exact values for $n\leq12$ and the following bounds.

\begin{theorem}
\emph{\cite{vb}} $\lfloor\frac{6n+9}{5}\rfloor\leq\gamma_{\mathrm{m}}^{\infty
}(P_{5}\boksie P_{n})\leq\lfloor\frac{4n+3}{3}\rfloor$.
\end{theorem}

We now compare the m-eternal domination number and the vertex cover number.
This may seem like an unusual pair of parameters to compare, but the
comparison turns out to be interesting.

\begin{theorem}
\label{gammavc}\label{tree-vc}

\begin{enumerate}
\item[\emph{(a)}] \emph{\cite{KMDM}\hspace{0.1in}}If $G$ is connected, then
$\gamma_{\mathrm{m}}^{\infty}(G)\leq2\tau(G)$.

\item[\emph{(b)}] \emph{\cite{KMDM}\hspace{0.1in}}If, in addition,
$\delta(G)\geq2$, then $\gamma_{\mathrm{m}}^{\infty}(G)\leq\tau(G)$.

\item[\emph{(c)}] \emph{\cite{KM6}\hspace{0.1in}}If, in addition to (a) and
(b), $G$ has girth seven or at least nine, then $\gamma_{\mathrm{m}}^{\infty
}(G)<\tau(G)$.

\item[\emph{(d)}] \emph{\cite{KM6}\hspace{0.1in}}For any nontrivial tree $T$,
$\tau(T)\leq\gamma_{\mathrm{m}}^{\infty}(T)\leq2\tau(T)$.
\end{enumerate}
\end{theorem}

It is not possible to relax the girth condition in Theorem \ref{tree-vc}(c) to
girth less than five. Examples of graphs with girth less than five for which
$\gamma_{\mathrm{m}}^{\infty}(G)=\tau(G)$ are given in \cite{KM6}. The problem
remains open for girths five, six, and eight, though it is believed that
$\gamma_{\mathrm{m}}^{\infty}(G)<\tau(G)$ for such graphs. The trees where the
bounds in Theorem \ref{tree-vc}(d) are sharp are characterized in \cite{KM6}.

A question stated in \cite{GHH} is whether there is any advantage in allowing
two guards to occupy the same vertex in the m-eternal domination problem.
There is no advantage allowing multiple guards to occupy a single vertex in
the \textquotedblleft one guard moves\textquotedblright\ model \cite{BCG2}.
The results stated up until now in this paper apply to the case when only one
guard is allowed to occupy each vertex. Finbow et al. have shown that there
exist graphs for which it is an advantage in the all-guards move model to
allow more than one guard on a vertex at a time \cite{FG}. We sketch a proof
of this, using $\gamma_{\mathrm{m}}^{*\infty}(G)$ to denote the number of
guards needed if more than one guard is allowed on a vertex at a time (and all
guards are allowed to move in response to an attack).

\begin{theorem}
There exists a graph $G$ such that $\gamma_{\mathrm{m}}^{*\infty}(G)=9$ and
$\gamma_{\mathrm{m}}^{\infty}(G)=10$.
\end{theorem}

\noindent\textbf{Proof}.\hspace{0.1in} (Sketch) Let $K_{4}-e$ be the graph
formed from $K_{4}$ by deleting an edge. Define a \textit{widget} to be the
graph formed by taking two $K_{4}-e$'s and combining one degree three vertex
from each into a single vertex (so a widget has seven vertices and two
vertices of degree three). Form graph $G_{5}$ by taking five widgets along
with an additional vertex $x$; add an edge between $x$ and the vertices of
degree three in each widget.

To see that $\gamma_{\mathrm{m}}^{\infty}(G)=10$, one can observe that there
must be two guards in a widget at some point in time: even though the
domination number of a widget is 1, because the degree two vertices in a
widget are independent and not adjacent to $x$ and the dominating vertex in a
widget is not adjacent to $x$, it follows that we need at least two guards in
each widget at all times there is not a guard on $x$. If there is a guard on
$x$, then at most one widget can contain one guard.

On the other hand, by maintaining at least one guard in each widget at all
times and two guards on $x$, which move in and out of attacked widgets, one
can see that nine guards suffice to protect $G_{5}$. This is done by moving
both guards from $x$ to the widget where an attack occurs (so that widget
contains three guards immediately after an attack) and two guards from the
previously attacked widget move to $x$, whilst the remaining guard in that
widget moves to the degree six vertex in the widget. \hfill$\blacksquare$

\medskip

The proof in \cite{FG} is more general than the sketch given above and shows
there are graphs where $\gamma_{\mathrm{m}}^{*\infty}(G)$ and $\gamma
_{\mathrm{m}}^{\infty}(G)$ can differ by any additive constant. It remains
open to prove whether or not for all graphs $G$ there is a constant $c > 1$
such that $c\gamma_{\mathrm{m}}^{*\infty}(G) \geq\gamma_{\mathrm{m}}^{\infty
}(G)$.

If any number of guards per vertex are allowed, then the bound in Theorem
\ref{princeupper} can be improved to $\lceil\frac{n}{2}\rceil-1$ when
$\delta(G)\geq2$ (with four small exceptions) \cite{prince2}. It is not known
whether their result holds if each vertex contains at most one guard. Under
these conditions Nordhaus-Gaddum results were also shown in \cite{prince2},
for example the following bound; they also characterize the graphs for which
equality holds.

\begin{theorem}
\label{princenord}\emph{\cite{prince2}} $\gamma_{\mathrm{m}}^{\infty
}(G)+\gamma_{m}^{\infty}(\overline{G})\leq n+1$.
\end{theorem}

\section{Eternal total domination}

\label{Sec_ETD}

Some results on eternal total domination are reviewed in this section. The
first result applies to the \textquotedblleft one-guard
moves\textquotedblright\ model.

\begin{theorem}
\emph{\cite{KM4}} For all graphs $G=(V,E)$ without isolated vertices,

\begin{enumerate}
\item[\emph{(a)}] $\gamma_{t}^{\infty}(G)>\gamma^{\infty}(G)$

\item[\emph{(b)}] $\gamma_{t}^{\infty}(G)\leq\gamma^{\infty}(G)+\gamma
(G)\leq2\gamma^{\infty}(G)\leq2\theta(G).$
\end{enumerate}
\end{theorem}

Klostermeyer and Mynhardt give a number of results on eternal total domination
\cite{KM4} in the all-guards move model, such as the following.

\begin{theorem}
\emph{\cite{KM4}} \label{allmove} For all graphs $G=(V,E)$ without isolated
vertices, $\gamma_{\mathrm{\operatorname{mt}}}^{\infty}(G)\leq2\gamma(G)$.
\end{theorem}

Results from \cite{GK2} focus on grid graphs and include the following.

\begin{theorem}
\emph{\cite{GK2}\hspace{0.1in}}

\begin{enumerate}
\item[\emph{(a)}] For any $n\geq3$, $\gamma_{\mathrm{\operatorname{mt}}%
}^{\infty}(P_{2}\boksie P_{n})=\lfloor\frac{2n}{3}\rfloor+2$.

\item[\emph{(b)}] For all $n\geq1$, $\gamma_{\operatorname{mt}}^{\infty}%
(P_{3}\boksie P_{n})=n+1$.

\item[\emph{(c)}] For any $n\geq1$, $\gamma_{\operatorname{mt}}^{\infty}%
(P_{4}\boksie P_{n})\leq\left\lfloor \frac{4n}{3}\right\rfloor +2.$
\end{enumerate}
\end{theorem}

Achieving good bounds for larger grids graphs seems quite difficult; by
\textquotedblleft`good\textquotedblright\ bounds we mean better than simply
partitioning the grid into disjoint, say $3\times n$, grids.

\section{Eternal vertex covering}

\label{Sec_EVC}

We emphasize that eternal vertex covering is non-trivial only for the
all-guards move model and thus our attention is limited to that model. Some
simple examples are as follows: $\tau_{\mathrm{m}}^{\infty}(C_{4})=2$,
$\tau_{\mathrm{m}}^{\infty}(C_{5})=3$ and $\tau_{\mathrm{m}}^{\infty}%
(P_{n})=2\tau(P_{n})$ if $n$ is odd \cite{KM5}. A fundamental bound is given next.

\begin{theorem}
\emph{\cite{KM5}}\label{alphabounds}\hspace{0.1in}For any nontrivial connected
graph $G$, $\tau(G)\leq\tau_{\mathrm{m}}^{\infty}(G)\leq2\tau(G)$.
\end{theorem}

Graphs satisfying the upper bound in Theorem \ref{alphabounds} are
characterized in \cite{KM5}. Some graphs where the lower bound is sharp are
described next.

\begin{proposition}
Each graph in the following classes satisfies $\tau_{\mathrm{m}}^{\infty
}(G)=\tau(G)$.

\begin{enumerate}
\item[\emph{(a)}] $K_{n}$

\item[\emph{(b)}] Petersen graph

\item[\emph{(c)}] $K_{m}\boksie K_{n}$

\item[\emph{(d)}] $C_{m}\boksie C_{n}$

\item[\emph{(e)}] Circulant graphs (to repel an attack along the edge $uv$,
move (say) the guard on $u$ to $v$ and move each other guard along its
incident edge that corresponds to $uv$ in the same orientation of the cycle).
\end{enumerate}
\end{proposition}

We next give some exact bounds for trees and grid graphs. Let $L$ denote the
number of leaves of a tree $T$.

\begin{theorem}
\emph{\cite{KM5}} \label{tree}For any nontrivial tree $T$, $\tau_{\mathrm{m}%
}^{\infty}(T)=|V-L|+1$.
\end{theorem}

\begin{theorem}
\label{grids}

\begin{enumerate}
\item[(a)] $\tau_{\mathrm{m}}^{\infty}(P_{1}\boksie P_{n})=n-1$.\vspace
{-0.08in}

\item[\emph{(b)}] If $n$ is even, then $\tau_{\mathrm{m}}^{\infty}%
(P_{n}\boksie P_{m})=\frac{nm}{2}=\tau(P_{n}\boksie P_{m})$.\vspace{-0.08in}

\item[\emph{(c)}] If $n,m>1$ are odd, $n\geq m$, then $\tau_{\mathrm{m}%
}^{\infty}(P_{n}\boksie P_{m})=\lceil\frac{nm}{2}\rceil=\tau(P_{n}\boksie
P_{m})+1$.
\end{enumerate}
\end{theorem}

We next compare $\tau_{m}^{\infty}$ with some of the other graph protection parameters.

\begin{theorem}
\emph{\cite{KGTN}\hspace{0.1in}}If $G$ is connected, then $\tau_{m}^{\infty
}(G)=\gamma(G)$ if and only if $G\in\{C_{4},K_{2}\}$.
\end{theorem}

\begin{theorem}
\emph{\cite{KM5}}\hspace{0.1in}If $G\neq C_{4}$ is a connected graph of order
$n\geq3$ with $\delta(G)\geq2$, then $\gamma_{\mathrm{m}}^{\infty}%
(G)<\tau_{\mathrm{m}}^{\infty}(G)$.
\end{theorem}

It seems a challenging problem to describe graphs with pendant vertices and
$\gamma_{\mathrm{m}}^{\infty}(G)=\tau_{\mathrm{m}}^{\infty}(G)$. Some examples
are given next. Part (a) of Proposition \ref{addnminus2} is from \cite{KM5}
and we thank Michael Fisher for pointing out the examples in the proof of part (b).

\begin{proposition}
\label{addnminus2} Let $G$ be a $2$-connected graph with $n$ vertices. Let
$G^{\prime}$ be a graph obtained from $G$ by attaching a pendant vertex to
each vertex of $G$ except the two vertices $u,v$.\label{here}

\begin{enumerate}
\item[\emph{(a)}] If $uv\in E$ then $\alpha_{\mathrm{m}}^{\infty}(G^{\prime
})=n$ and $\gamma_{\mathrm{m}}^{\infty}(G^{\prime})=n-1$.

\item[\emph{(b)}] If $uv\notin E$ then $\alpha_{\mathrm{m}}^{\infty}%
(G^{\prime})\geq n-1=\gamma_{\mathrm{m}}^{\infty}(G^{\prime})$.
\end{enumerate}
\end{proposition}

\noindent\textbf{Proof}.\hspace{0.1in}(a)$\hspace{0.1in}$Suppose we could
eternally defend the edges of $G^{\prime}$ with $n-1$ guards. Let $x\in
V(G)-\{u,v\}$ and let $y$ be the pendant vertex attached to $x$. We can force
guards onto both vertices $x,y$. Since each end-vertex is dominated, the edge
$uv$ is not protected. It is easy to see that the vertices of $G^{\prime}$ can
be protected by $n-1$ guards.

\noindent(b)\hspace{0.1in}Similar to (a), $n-2$ guards do not protect the
edges of $G^{\prime}$. To see that $n$ guards suffice to defend the edges,
initially place guards on the vertices of $G$ and then maintain at most one
guard on a pendant edge at any time. Letting $G=C_{5}$ is an example where
$\alpha_{\mathrm{m}}^{\infty}(G^{\prime})=n$ and $G=C_{4}$ is an example where
$\alpha_{\mathrm{m}}^{\infty}(G^{\prime})=n-1$. \hfill$\blacksquare$

\begin{proposition}
\emph{\cite{KM5}} \label{add-ones} Let $G$ be a $2$-connected graph with
$n\geq3$ vertices. Add one pendant vertex to $n-1$ vertices of $G$ and call
the resulting graph $G^{\prime}$. Then $\alpha_{\mathrm{m}}^{\infty}%
(G^{\prime})=\gamma_{\mathrm{m}}^{\infty}(G^{\prime})=n$.
\end{proposition}

It is an open question whether the condition of $G$ being 2-connected in
Proposition \ref{add-ones} can be replaced by minimum degree two.

An analog of realizable triples can be defined for edge protection. Results on
graphs $G$ having realizable triples ($\tau(G)$, $\tau_{\mathrm{m}}^{\infty
}(G)$, $\tau_{\mathrm{mt}}^{\infty}(G)$), where $\tau_{\mathrm{mt}}^{\infty
}(G)$ is the \emph{total eternal vertex cover}, are given in \cite{ABC2,
KGTN}. Any such realizable triple must satisfy the basic bound that for a
connected graph $G$ with more than two vertices, $\tau_{mt}^{\infty}%
(G)<2\tau_{m}^{\infty}(G)$ \cite{KGTN}. In \cite{ABC2} it is shown that
$\tau_{mt}^{\infty}(G)\leq\tau_{c}(G)+1\leq2\tau(G)$, where $\tau_{c}$ is the
size of a smallest connected vertex cover of $G$.

It is shown in \cite{gaspers2} that there exist graphs for which allowing
multiple guards to reside on a vertex at the same time reduces the number of
guards needed to defend the edges of the graph, in comparison to the eternal
vertex cover number. These authors leave obtaining good bounds on $k$ in the
following statement as an open problem:

\begin{center}
If one can defend any sequence of $k$ attacks on edges,

then one can defend any infinite sequence of attacks on edges.
\end{center}

\noindent Partial results on this question are given in \cite{ABC2}. For instance:

\begin{theorem}
\emph{\cite{ABC2}} If $T$ is a tree with $n-L$ guards, then there exists a
strategy to defend $V(T)$ attacks on the edges of $T$. That is, an adversary
can be forced to make $V(T)$ attacks before winning the eternal vertex cover game.
\end{theorem}

An alternate type of eternal vertex cover problem in which attacks are at
vertices while a vertex cover must be maintained at all times is explored in
\cite{KVC}.

\section{Other models}

\label{Sec_other} Eternal independent dominating sets were studied in
\cite{HM} and secure independent sets (analogous to secure dominating sets)
were studied in \cite{regan}.

\subsection{Eviction Model}

In the eviction model, each configuration $D_{i}$,\ $i\geq1$, of guards is
required to be a dominating set. An attack occurs at a vertex $r_{i}\in D_{i}$
such that there exists at least one $v\in N(r_{i})$ with $v\notin D_{i}$. The
next guard configuration $D_{i+1}$ is obtained from $D_{i}$ by moving the
guard from $r_{i}$ to a vertex $v\in N(r_{i})$, $v\notin D_{i}$ (i.e., this is
the \textquotedblleft one-guard moves\textquotedblright\ model). The size of a
smallest eternal dominating set in the eviction model for $G$ is denoted
$e^{\infty}(G)$. Simply put, attacks occur at vertices with guards and we must
move that guard to an unoccupied neighboring vertex. An attacked vertex is
required to have at least one neighboring vertex with no guard, otherwise
there would be no place for the guard to go.

This problem models a problem in computer networks where copies of a file are
stored throughout the network and files must sometimes be moved, or migrated,
due to maintenance at the server at which they are located. The goal is to
ensure a copy of the file is close to every vertex in the network. That is,
the locations of the files induce a dominating set at all times. The eviction
problem was introduced in \cite{KLM} and \textquotedblleft one-guard
moves\textquotedblright\ and \textquotedblleft all-guards
move\textquotedblright\ versions were defined. Most of the results in that
paper are for the one-guard moves model and we focus our attention to that
model here.

Some easy examples to illustrate the concept are $e^{\infty}(K_{1, m})=m$,
$e^{\infty}(C_{5})=2$, and $e^{\infty}(P_{5})=3$.

\begin{theorem}
\emph{\cite{KLM}\hspace{0.1in}}Let $G$ be a connected graph. Then $e^{\infty
}(G)\leq\theta(G)$.
\end{theorem}

\begin{theorem}
\emph{\cite{KLM}\hspace{0.1in}}Let $G$ be a bipartite graph. Then $e^{\infty
}(G)=\alpha(G)$.
\end{theorem}

Unlike in the traditional eternal domination model, there are graphs $G$ for
which $e^{\infty}(G)<\alpha(G)$: take a copy of $K_{2}$ and a large
independent set $I$ and join every vertex of the $K_{2}$ to every vertex of
$I$. This graph has $e^{\infty}(G)=1$.

\begin{theorem}
\emph{\cite{KLM}\hspace{0.1in}}There exists a graph $G$ such that $e^{\infty
}(G)>\alpha(G)$. In fact, for $k \geq3$, $e^{\infty}(C_{2k+1}) = k+1$.
\end{theorem}

\begin{theorem}
\cite{KLM} Let $G$ be a graph with $\alpha(G) = 2$. If $G$ has two dominating
vertices, then $e^{\infty}(G) = 1$. Otherwise, $e^{\infty}(G) = 2$.
\label{ind-bound2}
\end{theorem}

\noindent\textit{Proof:} If $G$ has dominating vertices $x$ and $y$, then a
single guard can relocate back and forth between them and maintain a
dominating set.

Finally, suppose $G$ has at most one dominating vertex. Then $G$ is the
complement of a triangle-free graph with at most one isolated vertex.
Initially locate the guards on any dominating set of size two, say $\{u, v\}$.
Suppose the guard on $u$ is attacked. If $v$ has a non-neighbor $w \neq u$,
then whether or not $u$ and $v$ are adjacent, the guard at $u$ guard can
relocate to $v$ and the resulting configuration is a dominating set. If no
such vertex $w$ exists, the guard at $u$ can relocate to any vertex $z$ and
the resulting configuration of guards is a dominating set. $\Box$

\medskip

The following result is much more difficult to prove.

\begin{theorem}
\emph{\cite{KLM}} Let graph $G$ have $\alpha(G)=3$. Then $e^{\infty}(G)\leq5$.
\end{theorem}

It is not known whether or not $e^{\infty}(G)\leq\gamma^{\infty}(G)$ for all
graphs $G$ \cite{KLM}.

Much less is known about the eviction model when all guards are allowed to
move in response to an attack, though some elementary results are given in
\cite{KLM} and in \cite{KMe}.

The eviction model for eternal independent sets was studied in \cite{CK}.

\subsection{Eternal Connected Domination}

Let $\gamma_{c}^{\infty}(G)$ denote the size of a smallest \emph{eternal
connected dominating set }(\emph{ECDS}) in which the vertices containing
guards induce a connected graph. Denote the all-guards move version of this
parameter (the cardinality of a minimum m-eternal\emph{ connected dominating
set} ($m$-\emph{ECDS}) by $\gamma_{\mathrm{\operatorname{mc}}}^{\infty}(G)$.
The ordinary connected domination number of $G$ is denoted $\gamma_{c}(G)$
\cite{HHS}. Obviously, these parameters are only defined for connected graphs.
They were initially studied in \cite{KM4}.

\begin{theorem}
\emph{\cite{KM4}\hspace{0.1in}}If $G$ is connected and $\theta(G)\geq2$, then
$\gamma_{\mathrm{mc}}^{\infty}(G)\leq2\theta(G)-1$. This bound is sharp for
all $\theta\geq2$.
\end{theorem}

\begin{theorem}
\emph{\cite{KM4}\hspace{0.1in}}For all graphs $G=(V,E)$ without isolated vertices,

\begin{enumerate}
\item[\emph{(a)}] $\gamma_{c}^{\infty}(G)>\gamma^{\infty}(G)$

\item[\emph{(b)}] $\gamma_{c}^{\infty}(G)\leq\gamma^{\infty}(G)+\gamma
(G)\leq2\gamma^{\infty}(G)\leq2\theta(G)$.
\end{enumerate}
\end{theorem}

Klostermeyer and Mynhardt also give a number of results on eternal connected
domination \cite{KM4} in the all-guards move model, such as the following bound.

\begin{theorem}
\emph{\cite{KM4}}\hspace{0.1in}For all graphs $G=(V,E)$ without isolated
vertices, $\gamma_{\mathrm{\operatorname{mc}}}^{\infty}(G)\leq2\gamma(G)$.
\end{theorem}

\subsection{Foolproof Eternal Domination}

In the definition of eternal domination, the decision of which guard to send
to defend an attack may require knowledge of the locations of future attacks.
The definition states \textquotedblleft there exists\textquotedblright\ a
guard to send to defend the attack such that all subsequent attacks can be
defended by the resulting guard configuration. It may be difficult in practice
to decide which guard to send to defend an attack.

Burger et al. \cite{BCG2} defined a \textquotedblleft foolproof" variation on
eternal domination in which the resulting configuration of guards must be able
to defend all subsequent attacks if a guard from \textbf{{any}} vertex
adjacent to the attacked vertex is sent to defend an attack at an unoccupied
vertex. That is, no matter which guard is sent, the resulting configuration
can defend all future attacks. They proved that $n-\delta(G)$ guards are
necessary and sufficient for all graphs $G$, where $\delta(G)$ is the minimum
vertex degree in the graph. To see this, note that any set of $n-\delta(G)$
vertices form a dominating set. On the other hand, if we have fewer than
$n-\delta(G)$ guards in $G$, then by a series of attacks, an adversary can
force the closed neighborhood of a vertex to contain no guards. For example,
consider $C_{6}$, and observe that $\gamma^{\infty}(C_{6})=3$. Now suppose we
could defend the graph with three guards in the foolproof model. Since any
neighboring guard can move to defend an attack, an adversary can force the
three guards to migrate to three consecutive vertices, thereby leaving a
vertex undominated.

The foolproof variety has been studied in the all-guards move model in
\cite{kmfool}. The problem is the same as the m-eternal dominating set problem
in that attacks are at (unoccupied) vertices and all guards can move in
response to an attack on a vertex $v$, but the attacker chooses which guard
moves to $v$. One can also imagine there being a victim of the attack at $v$
and allowing the victim to choose which guard to send to its defense. For
example, when a site is attacked, it may want to choose which of the nearby
defenders it calls in, perhaps because of particular expertise in defending
certain types of attacks. The size of a smallest m-eternal dominating set for
$G$ in the foolproof model is denoted $\rho_{m}^{\infty}(G)$.

\begin{proposition}
\emph{\cite{kmfool}\hspace{0.1in}}For any graph $G$, $\gamma_{\mathrm{m}%
}^{\infty}(G)\leq\rho_{m}^{\infty}(G)\leq\alpha_{m}^{\infty}(G)$.
\end{proposition}

\noindent\textbf{Proof}.\hspace{0.1in}The first inequality is obvious. For the
second inequality, observe that in the m-eternal vertex cover problem, when an
attack occurs on an edge with guards on either end, the two guards can swap
places and no other guards need to move; hence there is no net change in the
guard configuration. If there is only one guard incident to attacked edge
$uv$, that guard must move across the edge, say from $u$ to $v$, to defend the
attack. This is equivalent to the attacker choosing the guard to defend the
attack. Now rather than having attacks at edges, imagine the attack is at $v$
and the attacker chooses the guard at $u$ to defend it. It follows that
$\rho_{m}^{\infty}(G)\leq\alpha_{m}^{\infty}(G)$. \hfill$\blacksquare$

\begin{theorem}
\emph{\cite{kmfool}}\label{cliques}

\begin{enumerate}
\item[\emph{(a)}] If $G$ is a connected bipartite graph, then $\rho
_{m}^{\infty}(G)$.

\item[\emph{(b)}] For any graph $G$, $\rho_{m}^{\infty}(G)\leq2\theta(G).$
\end{enumerate}
\end{theorem}

It is not known if the bound in Theorem \ref{cliques} (b) is sharp. There does
exist a graph $G$ with $\rho_{m}^{\infty}(G) \geq\frac{3}{2} \theta(G)$
\cite{kmfool}.

\section{Complexity}

The complexity of deciding whether a given set of vertices is an eternal
dominating set, or another of the variations discussed, as well as the
complexity of determining the protection parameters themselves, are generally
difficult problems. The precise complexity remains unknown in most cases. For
example, it is unknown whether deciding whether a given set of vertices is an
eternal dominating set lies in the class PSPACE (though it is not too
difficult to see that it can be decided in exponential time, based on the
``configuration graph'' idea from \cite{KLM}). One problem in assessing in
which complexity class the eternal domination problem lies is to determine how
many attacks one must evaluate to determine whether a set of guards can defend
any infinite sequence of attacks in the graph. That is, is there a polynomial
function $f(n)$, where $n$ is the number of vertices in $G$, such that if one
can defend any sequence of $f(n)$ attacks, then one can defend any infinite
sequence of attacks? If there is no such polynomial function, then what bounds
can be placed on such a function?

We mention some results, besides the obvious cases like for perfect graphs.
From the results in \cite{KM2}, the $\mathrm{m}$-eternal domination number for
a tree can be computed in polynomial time. In addition, we can determine in
polynomial time whether each of these protection parameters is at most $k$ for
a fixed constant $k$, based the configuration graph method of \cite{KLM}. On a
related note, the parameterized complexity of the eternal vertex cover problem
was studied in \cite{gaspers}.

\section{Open problems}

\label{Sec_Open}

We present a number of conjectures and open problems on some of the models
discussed above.

\subsection{Eternal domination}

\begin{problem}
Study classes of graphs $G$ such that $(i)$\hspace{0.1in}$\gamma^{\infty
}(G)=\alpha(G)$, $(ii)$\hspace{0.1in}$\gamma^{\infty}(G)=\theta(G)$.
\end{problem}

As mentioned above, $\gamma^{\infty}(G)=\theta(G)$ if $G$ is series-parallel,
so it makes sense to pose the following question.

\begin{problem}
Is it true that $\gamma^{\infty}(G)=\theta(G)$ if $G$ is planar?
\end{problem}

\begin{problem}
Does there exist a constant $c$ such that $\gamma^{\infty}(G)\leq c\tau(G)$
for all graphs $G$?
\end{problem}

The following is motivated by an error discovered in \cite{KM2}, where it is
claimed that no such graph exists.

\begin{problem}
\label{counterexm} Does there exist a graph $G$ with $\gamma(G)=\gamma
^{\infty}(G)$ and $\gamma(G)<\theta(G)$?
\end{problem}

In \cite{KM8}, it was shown that (i) every triangle-free $G$ with
$\gamma(G)=\gamma^{\infty}(G)$ has $\gamma(G)=\theta(G)$ and (ii) every graph
$G$ with $\Delta(G) \leq3$ with $\gamma(G)=\gamma^{\infty}(G)$ has
$\gamma(G)=\theta(G)$.

It would also be of interest to determine if the graph with 11 vertices given
in \cite{GHH} having $\gamma^{\infty}<\theta$ is the smallest such graph.

\begin{problem}
\label{gamma-prob}\emph{\cite{KM2}\hspace{0.1in}}Characterize graphs $G$ with
$\gamma(G)=\gamma^{\infty}(G)=\theta(G)$.
\end{problem}

It is not hard to argue that any graph $G$ satisfying $\gamma(G)=\gamma
^{\infty}(G)<\theta(G)$ contains a triangle.

\begin{problem}

\begin{enumerate}
\item[\emph{(a)}] Describe classes of graphs with $\gamma^{\infty}%
/\alpha>\frac{3}{2}$.

\item[\emph{(b)}] Characterize realizable triples with $\gamma^{\infty}%
/\alpha>\frac{3}{2}$.
\end{enumerate}
\end{problem}

A Vizing-like question was asked in \cite{KM8}.

\begin{problem}
It is true for all graphs $G$ and $H$ that $\gamma^{\infty}(G\boksie H)\geq
\gamma^{\infty}(G)\ast\gamma^{\infty}(H)$?
\end{problem}

Interestingly, such a Vizing-like condition was shown not to hold for all
graphs $G$ in the all-guards-move model in \cite{KM8}.

\subsection{m-Eternal Domination}

\label{SecOpen_m-eternal}Recall the inequality chain $\gamma(G)\leq
\gamma_{\mathrm{m}}^{\infty}(G)\leq\alpha(G)\leq\gamma^{\infty}(G)\leq
\theta(G)$ from Section \ref{Sec_mED}.

\begin{problem}
Describe classes of graphs having $\gamma(G)=\gamma_{\mathrm{m}}^{\infty}(G)$,
$\gamma^{\infty}(G)=\gamma_{\mathrm{m}}^{\infty}(G)$, $\gamma_{\mathrm{m}%
}^{\infty}(G)=\tau(G)$, or $\gamma_{\mathrm{m}}^{\infty}(G)=\alpha(G)$.
\end{problem}

As shown in \cite{BSL}, there exist connected Cayley graphs, necessarily of
non-Abelian groups, whose m-eternal domination numbers exceed their domination
numbers by one. This implies that there exist disconnected Cayley graphs $G$
such that $\gamma_{\mathrm{m}}^{\infty}(G)-\gamma(G)$ is an arbitrary positive
integer. The picture for connected Cayley graphs is not so clear.

\begin{problem}
Does there exist a connected Cayley graph $G$ such that $\gamma_{\mathrm{m}%
}^{\infty}(G)>\gamma(G)+1$?
\end{problem}

\begin{problem}
Find conditions under which the bound $\gamma_{\mathrm{m}}^{\infty}%
(G)\leq\lceil\frac{n}{2}\rceil$ in Theorem \ref{princeupper} can be improved,
and conditions under which equality holds.
\end{problem}

\begin{problem}
Determine the value of $\gamma_{\mathrm{m}}^{\infty}(P_{n}\boksie P_{m})$. In
particular, is $\gamma_{\mathrm{m}}^{\infty}(P_{n}\boksie P_{n})\leq
\gamma(P_{n}\boksie P_{n})+c$, for some constant $c$? (The latter is
conjectured to be true by S. Finbow and W. Klostermeyer, personal communication).
\end{problem}

\begin{conjecture}
\emph{\cite{GK2}\hspace{0.1in}}If $\gamma_{\mathrm{m}}^{\infty}(P_{3}\boksie
P_{n})\leq r$, then $\gamma_{\mathrm{m}}^{\infty}(P_{3}\boksie P_{n+1})\leq
r+1$.
\end{conjecture}

\begin{conjecture}
\emph{\cite{GK2}\hspace{0.1in}For }$n>9$, $\gamma_{\mathrm{m}}^{\infty}%
(P_{3}\boksie P_{n})=1+\left\lceil \frac{4n}{5}\right\rceil $.
\end{conjecture}

The latter conjecture has been nearly resolved by Finbow et al. in Theorem
\ref{3-grid}, as discussed above.

\subsection{Eternal Vertex Cover}

\begin{problem}
\emph{\cite{KM5}\hspace{0.1in}}For which (bipartite) graphs is $\tau
_{\mathrm{m}}^{\infty}(G)=\tau(G)$?
\end{problem}

\begin{problem}
\emph{\cite{KM5}\hspace{0.1in}}Do all vertex transitive graphs $G$ satisfy
$\tau_{\mathrm{m}}^{\infty}(G)=\tau(G)$?
\end{problem}

\begin{conjecture}
\emph{\cite{KM5}\hspace{0.1in}}Let $G$ and $H$ be graphs such that
$\tau_{\mathrm{m}}^{\infty}(G)=\tau(G)$ and $\tau_{\mathrm{m}}^{\infty
}(H)=\tau(H)$. Then $\tau_{\mathrm{m}}^{\infty}(G\boksie H)=\tau(G\boksie H).$
\end{conjecture}

\begin{conjecture}
\emph{\cite{KM5}\hspace{0.1in}}Let $G=(V,E)$ be a connected graph with
subgraph $H$ such that $\delta(H)\geq2$ and $\delta(G[V-V(H)])\geq2$. Then
$\tau_{\mathrm{m}}^{\infty}(G)\geq\tau_{\mathrm{m}}^{\infty}(H)+\tau
_{\mathrm{m}}^{\infty}(G[V-V(H)])$.
\end{conjecture}

\begin{problem}
\emph{\cite{KM5}} Characterize graphs that are edge-critical for eternal
vertex covering.
\end{problem}

If $e\in E(\overline{G})$, then possibly $\tau_{\mathrm{m}}^{\infty}%
(G+e)>\tau_{\mathrm{m}}^{\infty}(G)$ (such as $G=C_{4}$) or possibly
$\tau_{\mathrm{m}}^{\infty}(G+e)<\tau_{\mathrm{m}}^{\infty}(G)$. An example of
the latter is to let $G+e$ be the $2\times4$ grid graph laid out in the usual
manner (this graph has eternal vertex cover number four) and choose $e$ be the
middle edge on the upper $P_{4}$.

In general, vertex and edge criticality has not been studied for any of the
eternal protection parameters.

\subsection{Other models}

We mention some open problems in some of the other models discussed.

\begin{conjecture}
\emph{\cite{KM4}\hspace{0.1in}}For all connected graphs $G$ with
$\Delta(G)<n-1$, $\gamma_{c}^{\infty}(G)>\theta(G)$.
\end{conjecture}

\begin{problem}
\emph{\cite{KLM}\hspace{0.1in}}Is $e^{\infty}(G)\leq\gamma^{\infty}(G)$ for
all graphs $G$?
\end{problem}

The analogous problem in the all-guards-move model of eviction (in which an
attacked vertex may be occupied by another guard after the attack) is also
open. In the model in all-guards-move model of eviction in which an attacked
vertex must remain unoccupied until the next attack, it is sometimes the case
that more guards are needed in the eviction model than the traditional
all-guards-move model: $\gamma_{\mathrm{m}}^{\infty}(K_{1,m})=2$ when $m
\geq2$, but this graph requires $m$ guards in the eviction model in which an
attacked vertex must remain unoccupied until the next attack.

\begin{conjecture}
\emph{\cite{kmfool}\hspace{0.1in}}For a graph $G=(V,E)$ with $n$ vertices and
no isolated vertices, $\rho_{\mathrm{\operatorname{m}}}^{\infty}(G)\leq
\lceil\frac{n}{2}\rceil$.
\end{conjecture}

It was shown in \cite{kmfool} that $\rho_{\mathrm{\operatorname{m}}}^{\infty
}(G)\leq\lceil\frac{5n}{6}\rceil$, for all graphs $G$.

\section*{Bibliography}

\subsection*{Roman domination}

\begin{enumerate}
\item M.~Adabi, E.~E.~Targhi, N.~Jafari Rad, M.~S.~Moradi, Properties of
independent Roman domination in graphs. \emph{Australas.~J.~Combin.
}\textbf{52} (2012), 11--18.

\item M.~Atapour, S.~M.~Sheikholeslami, A.~Khodkar, Trees whose Roman
domination subdivision number is 2. \emph{Util.~Math}.~\textbf{82} (2010), 227--240.

\item M.~Atapour, S.~M.~Sheikholeslami, A.~Khodkar, Roman domination
subdivision number of graphs. \emph{Aequationes Math}.~\textbf{78} (2009), 237--245.

\item B.~Chaluvaraju, V.~Chaitra, Roman domination in odd and even graphs.
\emph{South East Asian J.~Math.~Math.~Sci.}~\textbf{10} (2011), 97--104.

\item E.~Chambers, W.~Kinnersly, N.~Prince, D.~West, Extremal problems for
Roman domination. \emph{SIAM J.~Discrete Math}.~\textbf{23} (2009), 1575--1586.

\item M.~Chellali, N.~Jafari Rad, A note on the independent Roman domination
in unicyclic graphs. \emph{Opuscula Math.}~\textbf{32} (2012), 715--718.

\item E.~J.~Cockayne, P.~A.~Dreyer, S.~M.~Hedetniemi, S.~T.~Hedetniemi, Roman
domination in graphs. \emph{Discrete Math}.~\textbf{278 }(2004), 11--12.

\item N.~Dehgardi, S.~M.~Sheikholeslami, L.~Volkmann, On the Roman $k$-bondage
number of a graph. \emph{AKCE Int.~J.~Graphs Comb.}~\textbf{8} (2011), 169--180.

\item E.~Ebrahimi Targhi, N.~Jafari Rad, L.~Volkmann, Unique response Roman
domination in graphs. \emph{Discrete Appl.~Math}.~\textbf{159} (2011), 1110--1117.

\item O.~Favaron, H.~Karami, R.~Khoeilar, S.~M.~Sheikholeslami, On the Roman
domination number of a graph. \emph{Discrete Math}.~\textbf{309} (2009), 3447--3451.

\item H.~Fernau, Roman domination: a parameterized perspective.
\emph{Int.~J.~Comput.~Math}. \textbf{85} (2008), 25--38.

\item X.~Fu, Y.~Yang, B.~Jiang, Roman domination in regular graphs.
\emph{Discrete Math}.~\textbf{309} (2009), 1528--1537.

\item S.~Fujita, M.~Furuya, Difference between 2-rainbow domination and Roman
domination in graphs. \emph{Discrete Appl.~Math}.~\textbf{161} (2013), 806--812.

\item A.~Hansberg, N.~Jafari Rad, L.~Volkmann, Characterization of Roman
domination critical unicyclic graphs. \emph{Util.~Math.~}\textbf{86} (2011), 129--146.

\item A.~Hansberg, L.~Volkmann, Upper bounds on the $k$-domination number and
the $k$-Roman domination number. \emph{Discrete Appl.~Math}.~\textbf{157}
(2009), 1634--1639.

\item M.~A.~Henning, A characterization of Roman trees. \emph{Discussiones
Math.~Graph Theory} \textbf{22} (2002), 225--234.

\item M.~A.~Henning, Defending the Roman Empire from multiple attacks.
\emph{Discrete Math. }\textbf{271} (2003), 101--115.

\item N.~Jafari Rad, A note on Roman domination in graphs. \emph{Util.~Math}%
.~\textbf{83} (2010), 305--312.

\item N.~Jafari Rad, Chun-Hung Liu, Trees with strong equality between the
Roman domination number and the unique response Roman domination number.
\emph{Australas.~J.~Combin.~}\textbf{54} (2012), 133--140.

\item N.~Jafari Rad, S.~M.~Sheikholeslami, Roman reinforcement in graphs.
\emph{Bull.~Inst.~Combin.~Appl}.~\textbf{61} (2011), 81--90.

\item N.~Jafari Rad, L.~Volkmann, Roman domination dot-critical trees.
\emph{AKCE Int.~J. Graphs Comb}.~\textbf{8} (2011), 75--83.

\item N.~Jafari Rad, L.~Volkmann, On the Roman bondage number of planar
graphs. \emph{Graphs Combin}.~\textbf{27} (2011), 531--538.

\item N.~Jafari Rad, L.~Volkmann, Roman domination perfect graphs.
\emph{An.~\c{S}tiin\c{t}.~Univ.~"Ovidius\textquotedblright\ Constan\c{t}a
Ser.~Mat.}~\textbf{19} (2011), 167--174.

\item A.~P.~Kazemi, Roman domination and Mycieleki's structure in
graphs.\emph{ Ars Combin}.~\textbf{106} (2012), 277--287.

\item A.~Khodkar, B.~P.~Mobaraky, S.~M.~Sheikholeslami, Upper bounds for the
Roman domination subdivision number of a graph. \emph{AKCE Int.~J.~Graphs
Comb}.~\textbf{5} (2008), 7--14.

\item D.~Kuziak, Interpolation properties of domination parameters of a graph.
\emph{Australas.~J.~Combin}.~\textbf{46} (2010), 3--12.

\item S.~C.~Li, Z.~X.~Zhu, Roman domination numbers.~(Chinese)
\emph{J.~Huazhong Univ.~Sci. Technol.~Nat.~Sci}.~\textbf{34} (2006), 118--119.

\item M.~Liedloff, T.~Kloks, J.~Liu, S-L.~Peng, Efficient algorithms for Roman
domination on some classes of graphs. \emph{Discrete Appl.~Math}.~\textbf{156}
(2008), 3400--3415.

\item M.~Liedloff, T.~Kloks, J.~Liu, S-L.~Peng, Roman domination over some
graph classes. \emph{Graph-theoretic concepts in computer science,} 103--114,
Lecture Notes in Comput.~Sci., \textbf{3787}, Springer, Berlin, 2005.

\item Chun-Hung Liu, G.~J.~Chang, Roman domination on 2-connected graphs.
\emph{SIAM J.~Discrete Math.}~\textbf{26} (2012), 193--205.

\item C.-H.~Liu, G.~J.~Chang, Upper bounds on Roman domination numbers of
graphs. \emph{Discrete Math.~}\textbf{312} (2012), 1386--1391.

\item B.~P.~Mobaraky, S.~M.~Sheikholeslami, Bounds on Roman domination numbers
of graphs. \emph{Mat.~Vesnik} \textbf{60} (2008), 247--253.

\item A.~Pagourtzis, P.~Penna, K.~Schlude, K.~Steinh\"{o}ofel, D.~Taylor,
P.~Widmayer, Server placements, Roman domination and other dominating set
variants. \emph{Proceedings 2002 IFIP Conference on Theoretical Computer
Science}, Montreal, Canada, 280--291.

\item P.~Pavli\v{c}, J.~\v{Z}erovnik, Roman domination number of the Cartesian
products of paths and cycles. \emph{Electron.~J.~Combin.} \textbf{19} (2012),
no. 3, Paper 19, 37 pp.

\item P.~R.~L.~Pushpam, T.~N.~M.~M.~Mai, On efficiently Roman dominatable
graphs. \emph{J.~Combin.~Math.~Combin.~Comput}.~\textbf{67} (2008), 49--58.

\item P.~R.~L.~Pushpam, T.~N.~M.~M.~Mai, Edge Roman domination in graphs.
\emph{J.~Combin.~Math.~Combin.~Comput}.~\textbf{69} (2009), 175--182.

\item R.~R.~Rubalcaba, P.~J.~Slater, Roman dominating influence parameters.
\emph{Discrete Math}.~\textbf{307} (2007), 3194--3200.

\item R.~R.~Rubalcaba, M.~Walsh, Fractional Roman domination. Proceedings of
the Thirty-Eighth Southeastern International Conference on Combinatorics,
Graph Theory and Computing. \emph{Congr.~Numer}.~\textbf{187} (2007), 8--20.

\item W.~Shang, X.~Wang, X.~Hu, Roman domination and its variants in unit disk
graphs. \emph{Discrete Math.~Algorithms Appl}.~\textbf{2} (2010), 99--105.

\item S.~M.~Sheikholeslami, L.~Volkmann, The Roman domination number of a
digraph. \emph{Acta Univ.~Apulensis Math.~Inform.}~\textbf{27} (2011), 77--86.

\item X.~Song, W.~Shang, Roman domination in a tree. \emph{Ars Combin}%
.~\textbf{98} (2011), 73--82.

\item X.~Song, X.~Wang, Roman domination number and domination number of a
tree. \emph{Chinese Quart.~J.~Math}.~\textbf{21} (2006), 358--367.

\item T.~K.~\v{S}umenjak, P.~Pavli\v{c}, On the Roman domination in the
lexicographic product of graphs. \emph{Discrete Appl.~Math.}~\textbf{160}
(2012), 2030--2036.

\item E.~E.~Targhi, N.~Jafari Rad, C.~M.~Mynhardt, Y.~Wu, Bounds for
independent Roman domination in graphs.
\emph{J.~Combin.~Math.~Combin.~Comput.}~\textbf{80} (2012), 351--365.

\item Y.~Wu, H.~Xing, Note on 2-rainbow domination and Roman domination in
graphs. \emph{Appl.~Math.~Lett.}~\textbf{23} (2010), 706--709.

\item H-M.~Xing, X.~Chen, X-G.~Chen, A note on Roman domination in graphs.
\emph{Discrete Math}.~\textbf{306} (2006), 3338--3340.

\item V.~Zverovich, A.~Poghosyan, On Roman, global and restrained domination
in graphs. \emph{Graphs Combin.} \textbf{27} (2011), 755--768.
\end{enumerate}

\subsection*{Secure domination}

\begin{enumerate}
\item S.~Benecke, E.~J.~Cockayne, C.~M.~Mynhardt, Secure total domination in
graphs. \emph{Util.~Math.~}\textbf{74} (2007), 247--259.

\item S.~Benecke, P.~J.~P.~Grobler, J.~H.~van Vuuren, Protection of complete
multipartite graphs. \emph{Util.~Math}.~\textbf{71} (2006), 161--168.

\item A.~P.~Burger, E.~J.~Cockayne, W.~R.~Gr\"{u}ndlingh, C.~M.~Mynhardt,
J.~H.~van Vuuren, W.~Winterbach, Finite order domination in graphs.
\emph{J.~Combin.~Math.~Combin.~Comput.~}\textbf{49} (2004), 159--175.

\item A.~P.~Burger, A.~P.~de Villiers, J.~H.~ van Vuuren, Two algorithms for
secure graph domination, \emph{J.~Combin.~Math.~Combin.~Comput.~}\textbf{85}
(2013), 321--339.

\item A.~P.~Burger, A.~P.~de Villiers, J.~H.~ van Vuuren, A linear algorithm
for secure domination in trees, \emph{Discrete Appl.~Math.} \textbf{171}
(2014), 15--27.

\item A.~P.~Burger, A.~P.~de Villiers, J.~H.~ van Vuuren, The cost of edge
failure with respect to secure graph domination, \emph{Util.~Math}.
\textbf{95} (2014), 329--339.

\item A.~P.~Burger, A.~P.~de Villiers, J.~H.~ van Vuuren, Edge stability in
secure graph domination. \emph{Discrete Mathematics and Theoretical Computer
Science} \textbf{17 }(2015), 103--122.

\item A.~P.~Burger, A.~P.~de Villiers, J.~H.~ van Vuuren, On minimum secure
dominating sets of graphs, \emph{Quaestiones Mathematicae}, to appear.

\item A.~P.~Burger, A.~P.~de Villiers, J.~H.~ van Vuuren, Edge criticality in
secure graph domination. Manuscript.

\item E.~J.~Cockayne, Irredundance, secure domination and maximum degree in
trees. \emph{Discrete Math}.~\textbf{307 }(2007), 12--17.

\item E.~J.~Cockayne, O.~Favaron, C.~M.~Mynhardt, Secure domination, weak
Roman domination and forbidden subgraphs. \emph{Bull.~Inst.~Combin.~Appl}%
.~\textbf{39} (2003), 87--100.

\item P.~J.~P.~Grobler, C.~M.~Mynhardt, Secure domination critical graphs.
\emph{Discrete Math}.~\textbf{309} (2009), 5820--5827.

\item W.~F.~Klostermeyer, C.~M.~Mynhardt, Secure domination and secure total
domination in graphs. \emph{Discuss.~Math.~Graph Theory} \textbf{28} (2008), 267--284.

\item C.~M.~Mynhardt, H.~C.~Swart, E.~Ungerer, Excellent trees and secure
domination. \emph{Util.~Math}.~\textbf{67} (2005), 255--267.
\end{enumerate}

\subsection*{Weak Roman domination}

\begin{enumerate}
\item E.J.~Cockayne, P.J.P.~Grobler, W.R.~Gr\"{u}ndlingh, J.~Munganga,
J.H.~van Vuuren, Protection of a graph. \emph{Utilitas Math.}~\textbf{67
}(2005), 19--32.

\item M.A.~Henning, S.T.~Hedetniemi, Defending the Roman Empire -- A new
strategy, \emph{Discrete Math.}~\textbf{266} (2003), 239--251.

\item P.~R.~L.~Pushpam, T.~N.~M.~M.~Mai, Weak Roman domination in graphs.
\emph{Discuss.~Math.~Graph Theory} \textbf{31} (2011), 115--128.
\end{enumerate}


\begin{thebibliography}{99}                                                                                               %


\bibitem {ABB}M.~Anderson, C.~Barrientos, R.~Brigham, J.~Carrington,
R.~Vitray, J.~Yellen, Maximum demand graphs for eternal security.
\emph{J}.~\emph{Combin.~Math.~Combin.~Comput.} \textbf{61} (2007), 111--128.

\bibitem {ABC2}M.~Anderson, R.~Brigham, J.~Carrington, R.~Dutton R.~Vitray,
J.~Yellen, Mortal and eternal vertex covers. Manuscript (2012).

\bibitem {ABC}M.~Anderson, R.~Brigham, J.~Carrington, R.~Dutton R.~Vitray,
J.~Yellen, Graphs simultaneously achieving three vertex cover numbers.
Manuscript (2012).

\bibitem {AF}J.~Arquilla, H.~Fredricksen, \textquotedblleft
Graphing\textquotedblright\ an optimal grand strategy. \emph{Military
Operations Research} \textbf{1}(3) (1995), 3--17.

\bibitem {BFM}I.~Beaton, S.~Finbow, J.~A.~MacDonald, Eternal domination of
grids. \emph{J}.~\emph{Combin.~Math.~Combin.~Comput.} \textbf{85} (2013), 33--48.

\bibitem {BFM2}I.~Beaton, S.~Finbow, J.~A.~MacDonald, Eternal domination of
grids. Manuscript (2013).

\bibitem {BC}B.~Bollobas, E.J.~Cockayne, Graph theoretic parameters concerning
domination, independence, irredundance. \emph{J. Graph Theory} \textbf{3}
(1979), 241--250.

\bibitem {BSL}A.~Braga, C.~C.~de Souza, O.~Lee, A note on the paper
\textquotedblleft Eternal security in graphs\textquotedblright\ by Goddard,
Hedetniemi, and Hedetniemi (2005). Manuscript (2014).

\bibitem {Braga2}A.~Braga, C.~C.~de Souza, O.~Lee, The Eternal Dominating Set
problem for proper interval graph, Manuscript (2015).

\bibitem {BCG2}A.~P.~Burger, E.~J.~Cockayne, W.~R.~Gr\"{u}ndlingh,
C.~M.~Mynhardt, J.~H.~van Vuuren, W.~Winterbach, Infinite order domination in
graphs. \emph{J}.~\emph{Combin.~Math.~Combin.~Comput.}~\textbf{50} (2004), 179--194.

\bibitem {CK}Y. Caro, W. Klostermeyer, Eternal Independent Sets in Graphs.
Manuscript (2015).

\bibitem {prince2}E.~Chambers, W.~Kinnersly, N.~Prince, Mobile eternal
security in graphs. Manuscript (2006).

\bibitem {CL}G.~Chartrand, L.~Lesniak, \emph{Graphs Digraphs}. Fourth Edition,
Chapman \& Hall, London, 2005.

\bibitem {ckpv}M.~Chrobak, H.~Karloff, T.~Payne, S.~Vishwanathan, New results
on server problems. \emph{SIAM J.~Discrete Math.}~\textbf{{4}} (1991), 172--181.

\bibitem {Protection}E.~J.~Cockayne, P.~J.~P.~Grobler, W.~R.~Grundlingh,
J.~Munganga, J.~H.~van Vuuren, Protection of a graph. \emph{Utilitas
Math}.~\textbf{67 }(2005), 19--32.

\bibitem {FG}S.~Finbow, S.~Gaspers, M.-E.~Messinger, P.~Ottoway, Eternal
domination. Manuscript (2010).

\bibitem {FMV}S.~Finbow, M.-E.~Messinger, M. van Bommel, Eternal domination in
$3 \times n$ grids.\emph{Australas.~J.~Combin.~}\textbf{61} (2015), 156--174.

\bibitem {gaspers}F.~Fomin, S.~Gaspers, P.~Golovach, D.~Kratsch, S.~Saurabh,
Parameterized algorithm for eternal vertex cover.
\emph{Inform.~Process.~Lett.}~\textbf{110} (2010), 702--706.

\bibitem {gaspers2}F.~Fomin, S.~Gaspers, P.~Golovach, D.~Kratsch, S.~Saurabh,
Eternal vertex cover. Manuscript (2010).

\bibitem {GHH}W.~Goddard, S.~M.~Hedetniemi, S.~T.~Hedetniemi, Eternal security
in graphs. \emph{J. Combin.~Math.~Combin.~Comput.}~\textbf{52} (2005), 169--180.

\bibitem {GK}J.~Goldwasser, W.~F.~Klostermeyer, Tight bounds for eternal
dominating sets in graphs. \emph{Discrete Math.~}\textbf{308} (2008), 2589--2593.

\bibitem {GK2}J.~Goldwasser, W.~F.~Klostermeyer, C.~M.~Mynhardt, Eternal
protection in grid graphs. \emph{Utilitas Math.}, ~\textbf{91} (2013), 47--64.

\bibitem {HM}B. Hartnell, C.~M.~Mynhardt, Independent Protection in Graphs.
Manuscript (2015).

\bibitem {HH2003}M.~A.~Henning, S.~T.~Hedetniemi, Defending the Roman Empire
-- a new strategy. \emph{Discrete Math}.~\textbf{266} (2003), 239--251.

\bibitem {HHS}T.~W.~Haynes, S.~T.~Hedetniemi, P.~J.~Slater, \emph{Fundamentals
of Domination in Graphs}. Marcel Dekker, New York, 1998.

\bibitem {KVC}W.~Klostermeyer, An eternal vertex cover problem. \emph{J.
Combin.~Math.~Combin.~Comput.}~\textbf{85} (2013), 79--95.

\bibitem {KGTN}W.~Klostermeyer, Some Questions on Graph Protection.
\emph{Graph Theory Notes of New York} \textbf{57} (2010), 29--33.

\bibitem {KLM}W.~Klostermeyer, M.~Lawrence, G.~MacGillivray, Dynamic
Dominating Sets: the Eviction Model for Eternal Domination. Manuscript (2014).

\bibitem {klos1}W.~F.~Klostermeyer, G.~MacGillivray, Eternally secure sets,
independence sets, and Cliques. \emph{AKCE Int.~J.~Graphs Comb}.\emph{~}%
\textbf{{2}} (2005), 119--122.

\bibitem {KM}W.~F.~Klostermeyer, G.~MacGillivray, Eternal security in graphs
of fixed independence number. \emph{J}.~\emph{Combin.~Math.~Combin.~Comput.}%
~\textbf{63} (2007), 97--101.

\bibitem {KM2}W.~F.~Klostermeyer, G.~MacGillivray, Eternal dominating sets in
graphs. \emph{J}.~\emph{Combin.~Math.~Combin.~Comput.} \textbf{68} (2009), 97--111.

\bibitem {kmfool}W.F.~Klostermeyer, G.~MacGillivray, Foolproof eternal
domination in the all-guards move model. \emph{Math Slovaca} \textbf{62}
(2012), 595--610.

\bibitem {KM-trees}W~.F.~Klostermeyer, G.~MacGillivray, Eternal domination in
trees, to appear in \textit{Journal of Combinatorial Mathematics and
Combinatorial Computing} (2014).

\bibitem {KM5}W.~F.~Klostermeyer, C.~M.~Mynhardt, Edge protection in graphs.
\emph{Australas.~J.~Combin.}~\textbf{45} (2009), 235--250.

\bibitem {KM4}W.~F.~Klostermeyer, C.~M.~Mynhardt, Eternal total domination in
graphs. \emph{Ars Combin.}~\textbf{68} (2012), 473--492.

\bibitem {KMDM}W.F.~Klostermeyer, C.M.~Mynhardt, Graphs with equal eternal
vertex cover and eternal domination numbers. \emph{Discrete Math.~}%
\textbf{311} (2011),1371--1379.

\bibitem {KM6}W.~F.~Klostermeyer, C.~M.~Mynhardt, Vertex covers and eternal
dominating sets. \emph{Discrete Appl.~Math.} \textbf{160} (2012), pp. 1183--1190.

\bibitem {KMe}W.~F.~Klostermeyer, C.~M.~Mynhardt, A dynamic domination problem
in trees. To appear in \emph{Transactions on Combinatorics} (2015).

\bibitem {KM8}W.~F.~Klostermeyer, C.~M.~Mynhardt, Domination, Eternal
Domination, and Clique Covering. To appear in \emph{Discuss.~Math.~Graph
Theory} (2015).

\bibitem {kout}E.~Koutsoupias, C.~Papadimitriou, On the $k$-server conjecture.
\emph{J.~ACM} \textbf{{42}} (1995), 971--983.

\bibitem {mms}M.~Manasse, L.~McGeoch, D.~Sleator, Competitive algorithms for
server problems. \emph{J.~Algorithms} \textbf{{11}} (1990), 208--230.

\bibitem {ms}L.~McGeoch, D.~Sleator, A strongly competitive randomized paging
algorithm. \emph{Algorithmica} \textbf{{6}} (1991), 816--825.

\bibitem {regan}F.~Regan, \emph{Dynamic variants of domination and
independence in graphs}, graduate thesis, Rheinischen Friedrich-Wilhlems
University, Bonn, 2007.

\bibitem {ReVelle}C.~S.~ReVelle, Can you protect the Roman Empire? \emph{Johns
Hopkins Magazine} \textbf{50}(2), April 1997.

\bibitem {ReVelleRos}C.~S.~ReVelle, K.~E.~Rosing, Defendens Imperium Romanum:
A classical problem in military strategy. \emph{Amer.~Math.~Monthly}
\textbf{107}, Aug.~-- Sept.~2000, 585--594.

\bibitem {Stewart}I.~Stewart, Defend the Roman Empire! \emph{Scientific
American,} December 1999, 136--138.

\bibitem {vb}C.M. van Bommel, M.F. van Bommel, Eternal Domination Numbers of
$5 \times n$ Grid Graphs, to appear in \textit{Journal of Combinatorial
Mathematics and Combinatorial Computing} (2015).

\bibitem {West}D.~B.~West, \emph{Introduction to Graph Theory}, Prentice-Hall, 1996.
\end{thebibliography}
\end{document}